\documentclass[article,moreauthors,pdftex]{mdpi} 
\usepackage{graphicx}
\usepackage{amsmath}
\usepackage{amssymb}
\usepackage{mathtools}
\usepackage{microtype}
\usepackage{bm}
\usepackage[utf8]{inputenc}
\allowdisplaybreaks
\preto{\abstractkeywords}{\nolinenumbers}
\newcommand{\ds}{\displaystyle }

\Title{Biased continuous-time random walks with Mittag-Leffler jumps}

\Author{Thomas M. Michelitsch$^{1}$, Federico Polito$^{2}$ and Alejandro P. Riascos$^{3}$}
\AuthorNames{Thomas M. Michelitsch, Federico Polito and  Alejandro P. Riascos}

\address{%
$^{1}$ \quad Sorbonne Universit\'e, Institut Jean le Rond d’Alembert, 
CNRS UMR 7190, 4 place Jussieu, 75252 Paris cedex 05, France; michel@lmm.jussieu.fr\\
$^{2}$ \quad Department of Mathematics ``Giuseppe Peano'', University of Torino, Italy; federico.polito@unito.it\\
$^{3}$ \quad Instituto de F\'isica, Universidad Nacional Aut\'onoma de M\'exico, Apartado Postal 20-364, 01000 Ciudad de M\'exico, M\'exico; aperezr@fisica.unam.mx}

\corres{Correspondence: michel@lmm.jussieu.fr}

\begin{document}

\maketitle

\begin{abstract}
\\
We construct admissible circulant
Laplacian matrix functions as generators for strictly increasing random walks on the integer line.
These Laplacian matrix functions refer to a certain class of Bernstein functions. 
The approach has connections with biased walks on digraphs.
Within this framework, we introduce a space-time generalization of the Poisson process as a strictly increasing walk with discrete 
Mittag-Leffler jumps time-changed with an independent (continuous-time) fractional Poisson process. 
We call this process `{\it space-time Mittag-Leffler process}'. We derive explicit formulae for the state probabilities which solve a Cauchy problem with a Kolmogorov-Feller (forward) difference-differential equation of general fractional type.
We analyze a 
``well-scaled'' diffusion limit and obtain a Cauchy problem with a space-time convolution equation involving Mittag-Leffler densities.
We deduce in this limit the `state density kernel' solving this Cauchy problem.
It turns out that the diffusion limit exhibits connections to Prabhakar general fractional calculus.
We also analyze in this way a generalization of the space-time Mittag-Leffler process.
The approach of constructing good Laplacian generator functions has a large potential in applications of space-time generalizations of the Poisson process and in the field of continuous-time random walks on digraphs.
\\[1mm]
\end{abstract}
\section{INTRODUCTION}
\label{Intro}
During the last three decades a vast interest in fractional calculus has emerged due to its success to describe so called `anomalous phenomena' such as
anomalous transport and diffusion \cite{KutnerMasoliver1990,GorenfloMainardi2006,Gorenflo2010, MetzlerKlafter2000,MetzlerKlafter2004,Zaslavsky2002,SaichevZaslavsky1997}, anomalous relaxation in dielectrics \cite{Capelas-Vaz2011}, creep models \cite{MainardiSpada2011},
and various complex phenomena \cite{GorenfloMainardi-complex-2006}, just to quote a few examples. Many of these models can be traced back to the Montroll-Weiss picture of continuous-time random walk (CTRW) \cite{MontrollWeiss1965}. Whereas classical CTRW models are random walks subordinated to an independent Poisson process, their fractional generalizations lead to
fat-tailed waiting-time Mittag-Leffler type densities with non-Markovian long memory behavior governed by evolution equations
of fractional or generalized fractional types
\cite{HilferAnton1995,MetzlerKlafter2000,MeerschaertEtal2011,Zaslavsky2002,KutnerMasoliver1990}.
For fundamental discussions on fractional calculus we invite the reader to consult the references
\cite{Ortiguera2015,Tarasov2018,Giusti2018a,HilferLuchko2019}.
In the meantime, several generalizations of fractional calculus have been proposed. One of the
most pertinent generalizations seems to be the so called Prabhakar generalization 
\cite{GaraGorenfloPolitoTomovski2014,MainardiGarrappa2015,GiustiPolitoMainardi-etal2020}. This generalization involves the Prabhakar function in place of the Mittag-Leffler function and was first introduced by Prabhakar \cite{Prabhakar1971}.
In the context of continuous-time renewal processes a Prabhakar generalization of the fractional Poisson process was introduced by Cahoy and Polito \cite{PolitoCahoy2013}, applied to stochastic motions on undirected graphs by Michelitsch and Riascos 
\cite{TMM-APR-PhysicaA2020,MichelitschRiascos2020,MiRia2020}, and discrete-time versions were developed in our recent paper \cite{MichelitschPolitoRiascos2020}.
\\[1mm]
The present paper is mainly concerned with space-time generalizations of the Poisson process
which can be traced back to special kinds of biased walks, namely {\it strictly increasing walks} on the integer line $\mathbb{N}_0$ subordinated to independent renewal processes. Such processes turn out to be governed by Cauchy problems with space-time difference-differential equations of general fractional type (generalized Kolmogorov-Feller or Kolmogorov-forward equation), 
\begin{equation}
\label{Cauchyproblemgen}
\begin{array}{clc}
 \ds 
{\cal D}_tp_n(t) & = \ds  -\frac{\xi}{g(1)}g(1-{\hat T}_{-1})p_n(t) &   \\ \\
\ds 
p_{n}(t)|_{t=0} &= \ds \delta_{n0} &
\end{array}\hspace{1cm} \xi>0,  \hspace{0.5cm} t \geq 0 ,\hspace{0.5cm} n \in \mathbb{N}_0,
\end{equation}
for the state probabilities $\{p_n(t)\}$ ($n=0,1,2,.. \in \mathbb{N}_0$) of the process. The state probabilities $p_n(t)={\mathbb P}({\cal N}(t)=n)$, $n \in \mathbb{N}_0$, $t \in \mathbb{R}^{+}$ indicate the probabilities to find a walker (with the indicated initial condition) at time $t$ in state ${\cal N}(t) = n$ in a strictly increasing walk 
${\cal N}(t)=\sum_{j=1}^{N(t)}Z_j \in \mathbb{N}_0$ with IID strictly positive integer jumps $Z_j\in \mathbb{N}$, where $N(t)\in \mathbb{N}_0$ denotes the number of arrivals up to time $t$ in a renewal process independent of the jumps $Z_j$. The state probability distribution is normalized on the state space $\sum_{n=0}^{\infty}p_n(t)=1 ,\,\, \forall t\geq 0$. On the left-hand side of \eqref{Cauchyproblemgen}, ${\cal D}_t$ stands for a general fractional derivative which is related to the mentioned renewal process (see \cite{Kochubei2011} for an outline of the framework of general fractional calculus). ${\hat T}_{-1}$ indicates the backward shift operator where shift operators ${\hat T}_{m}$ ($m\in \mathbb{Z}$) act on the {\it state space} such that ${\hat T}_{m}p_{n}(t)=p_{n+m}(t)$ (we denote
with $1$ the unity operator ${\hat T}_{0}=1$). The `Laplacian operator function' $g(1-{\hat T}_{-1})$ represents
the {\it generator} of the process. The stochastic process governed by Cauchy problem (\ref{Cauchyproblemgen}) is a space-time generalization of the Poisson process. Indeed for
$g(1-{\hat T}_{-1})=1-{\hat T}_{-1}$ (with 
$(1-{\hat T}_{-1})p_n(t)=p_n(t)-p_{n-1}(t)$) and ${\cal D}_t = \frac{d}{dt}$, Eq. (\ref{Cauchyproblemgen}) recovers 
the classical {\it Kolmogorov-Feller (forward) equation} of the homogeneous Poisson process. One of the goals of the present paper 
is to elaborate a general approach to construct non-trivial 
generalizations by means of ``good Laplacian operator functions'' $g(1-{\hat T}_{-1})$.
\\[1mm]
We will show that under suitable ``well-scaled'' continuous-space limit conditions from the integer-state space to a continuous one (i.e.\ from $\{0,1,\ldots\} =\mathbb{N}_0$ to $[0,\infty) = \mathbb{R}^{+}$) the discrete Cauchy problem (\ref{Cauchyproblemgen}) turns into a Cauchy problem with a diffusion equation of the following general structure:
\begin{equation}
\label{general}
\begin{array}{clc}
\ds {\cal D}_t{\cal P}(x,t)& = \ds  -\xi {\cal P}(x,t) + \xi \int_{0}^x{\cal W}(x-\tau){\cal P}(\tau,t){\rm d}\tau ,\hspace{1cm} \xi >0,\hspace{0.5cm} x,t \geq 0 &\\ \\
\ds {\cal P}(x,t)|_{t=0} &=  \delta(x) &
\end{array}
\end{equation}
This equation allows solely jumps into the positive $x$-direction, i.e. `strictly increasing' walks.
The convolution on the right-hand side describes incoming jumps from $[0,x)$
to state $x$ whereas the term $-\xi {\cal P}(x,t) =-\xi {\cal P}(x,t) \int_x^{\infty} {\cal W}(\tau-x,t){\rm d}\tau$ 
accounts for outgoing 
jumps from $x$ to $(x,\infty)$. Here ${\cal P}(x,t)$ stands for the `state density', ${\cal W}(x)$ indicates the transition density kernel and $\delta(x)$ the Dirac's $\delta$-distribution. 
\\[1mm]
A further aim of the present paper is to highlight connections of these models 
with a certain class of biased random walks on directed graphs. 
This aim is especially motivated by the huge upswing of
network science which
has become nowadays an immense interdisciplinary field \cite{Newman2010,NohRieger2004} last but not least driven by rapidly developing applications in online (social) networks and search engines. 
There is already a vast amount of specialized literature on various aspects of the subject, see e.g. \cite{Hughes1995,Hughes1996,Newman2010,Mohar1991lsg,Mohar1997sal}.
For instance, models have been developed to describe L\'evy flight dynamics and random walks with long-range jumps on undirected graphs \cite{TMM-APR-ISTE2019,RiascosMateos2014}, 
random walks with stochastic resetting on graphs \cite{RiascosBoyerHerringerMateos2019} and models for long-range mobility in cities \cite{RiascosMateos2020}, among many others. 
\\[1mm]
The present paper is organized as follows. In Section \ref{digraph} we recall some basic properties of biased walks on {\it directed networks}. We focus on ergodic (strongly connected) finite directed graphs and the features of their Laplacian and transition matrices. Having 
recalled these basic features, we define in Section \ref{biasedwalks} necessary and sufficient ``good Laplacian properties'' that a physically admissible `good Laplacian matrix' must fulfill. Starting with a good Laplacian matrix allowing only `local' transitions to next neighbor nodes we use these Laplacian properties to construct non-trivial good Laplacian functions (i.e. those that conserve the ``good Laplacian properties''). They will constitute the family of admissible generators $g(1-{\hat T}_{-1})$ on the right-hand side of (\ref{Cauchyproblemgen}). 
It turns out that the family of good Laplacian functions refers to a certain class of Bernstein functions.
For undirected graphs this approach was developed recently \cite{RiascosMichelitsch-et-al2018,TMM-APR-ISTE2019} and has also been extended to a class of biased walks \cite{RiascosMichelitschPizarro2020}. 
\\[1mm]
In section \ref{DTRW} we recall, within this picture, a class of biased walks on the integer line allowing solely strictly positive integer jumps generated by the fractional Laplacian matrix of a directed line. This leads to the {\it Sibuya walk}, 
which is characterized by a fat-tailed jump distribution, as a proto-typical example of a discrete approximation of 
the stable subordinator (see e.g. \cite{MeerschaertEtal2011,GorenfloMainardi2013,PachonPolitoRicciuti2018}).
\\[1mm]
Section \ref{space-time-gens} is devoted to explore, by means of this approach, space-time generalizations of 
the Poisson process leading to Cauchy problems of the general structure (\ref{Cauchyproblemgen}) involving ``good Laplacian matrix functions'' of Section \ref{biasedwalks}.
Within this framework, we also discuss classical cases such as the {\it space-fractional Poisson process} and the {\it space-time fractional Poisson process}, introduced by Orsingher and Polito \cite{OrshingerPolito2012}.
\\[1mm]
In the main part of this paper, Section \ref{space-time-frac}, we apply this approach to construct the 
``{\it space-time Mittag-Leffler process}'' as a pertinent generalization of the Poisson process governed by an evolution equation of general type (\ref{Cauchyproblemgen}). This process is a strictly increasing walk on the integer line with (discrete) Mittag-Leffler jumps separated by the Mittag-Leffler waiting times of the fractional Poisson process.
We derive explicitly the {\it state probabilities} and the {\it `well-scaled' continuous-space limit (diffusion limit)} of this process. We obtain a biased (forward) diffusion equation of general space-time fractional type
referring to the class (\ref{general}),
which is solved by the continuous-space limit density kernel of the state-probabilities. This kernel is derived in explicit form involving so called Prabhakar kernels. In this way, we show connections with Prabhakar general fractional calculus. 
We also highlight the equivalence with the Montroll-Weiss CTRW picture.
\\[1mm]
Finally, as a byproduct, we construct in Section \ref{Prabhakargeneralisation} the space-fractional generalization of the space-time Mittag-Leffler process and derive the Cauchy problem governing the state probabilities. The latter as well as its diffusion limit are derived in explicit forms where again Prabhakar kernels and general fractional calculus come into play.

\section{BIASED WALKS ON DIRECTED GRAPHS}
\label{digraph}
In this section we recall some basic notions and properties of biased walks taking place on directed graphs (also referred to as digraphs). We consider first a finite directed graph consisting of $N$ nodes (states) which we denote by $j=1,\ldots N$.
In a directed graph the edges between nodes have in general direction so that a path
may exist allowing the walk $m \to n$ but the return path $n \to m$ does
not necessarily exist. 
If for all pairs of nodes finite paths (and hence return paths) via directed edges exist, then 
the digraph is said to be {\it strongly connected} and equivalently {\it ergodic} (see e.g.  \cite{HararyPalmer1973,RiascosMichelitschPizarro2020}).
In order to define random walks on digraphs we introduce the $N\times N$ Laplacian 
matrix $\mathbf{L}$ containing this topological information. The Laplacian matrix has the elements \cite{RiascosMichelitschPizarro2020}
\begin{equation}
\label{Laplaciandigraph}
L_{ij}=k_i^{(\mathrm{out})} \delta_{ij}-\Omega_{ij}
\end{equation}
where we use the synonymous notation $\delta_{pq}=\delta_{p,q}$ for the Kronecker symbol.
Matrix (\ref{Laplaciandigraph}) is a generalization of the Laplacian matrix for binary undirected networks \cite{Newman2010,Mohar1991lsg,Mohar1997sal} to include the possibility of weights in the connections and asymmetry in the flow along the edges. For these particular connections we have $\Omega_{ij}\neq \Omega_{ji}$ where we assume $\Omega_{ij} \geq 0$ is a non-negative matrix.
Further, by construction we assume $\Omega_{ii}=0$.
In (\ref{Laplaciandigraph}) is present the so called {\it out-degree} 
defined by \cite{RiascosMichelitschPizarro2020}
\begin{equation}
\label{oudegree}
k_i^{(\mathrm{out})} = \sum_{j=1}^N \Omega_{ij} > 0 
\end{equation}
constituting a measure of the number of nodes which can be reached in a single jump from node $i$. The out-degree
is assumed to be {\it strictly positive} meaning that we do not allow isolated (disconnected) nodes.
Then, we introduce the {\it one-step transition matrix} $\mathbf{W}$  \cite{NohRieger2004,TMM-APR-ISTE2019,Newman2010}
\begin{equation}
\label{one-jump}
W_{i\to j} = \delta_{ij}- {\cal L}_{ij} = 
\frac{\Omega_{ij}}{k_i^{(\mathrm{out})}}
\end{equation}
denoting the probability of the transition $i\to j$ in one jump.
Per construction (\ref{one-jump}) is 
row-stochastic, i.e. $0\leq W_{i\to j} \leq 1$ with $\sum_{j=1}^N W_{i\to j} =1$ and the transition matrix is such that $W_{i\to i}=0$, i.e.\ in each jump the walker has to move to a different node.
For our convenience and as an equivalent description we introduce here the (non-symmetric) auxiliary Laplacian matrix
\begin{equation}
\label{reduced}
{\cal L}_{ij} = \frac{L_{ij}}{k_i^{(\mathrm{out})}}
\end{equation}
having normalized degrees ${\cal L}_{ii}=1$.
We notice that dynamical processes in directed networks have a greater variety than in the undirected case. Since the Laplacian matrices (\ref{Laplaciandigraph}) and (\ref{reduced}) are not symmetric
they have in general complex eigenvalues which are considered more closely in the subsequent Section \ref{biasedwalks}.
We can now define a discrete-time Markov chain (Markovian random walk) on this graph governed by
the simple master equation (e.g.\ \cite{NohRieger2004,TMM-APR-ISTE2019})
\begin{equation}
\label{masterq}
P_{ij}(t+1) = \sum_{r=1}^NP_{ir}(t)W_{r\to j} , \hspace{1cm} P_{ij}(0) =\delta_{ij}, \hspace{1cm} t \in \mathbb{N}_0.
\end{equation}
We assume here that the walker at $t=0$ is sitting on the departure node $i$. In the `$t$-step transition matrix' the element $P_{ij}(t)$ indicates the probability of the transition $i\to j$ in $t$ jumps. Note that for digraphs $\mathbf{P}(t)$ is not symmetric and for the initial condition $\mathbf{P}(0) =\mathbf{1}$ we have
\begin{equation}
\label{recursion}
[\mathbf{P}(t)]_{ij} = [\mathbf{W}^t]_{ij} ,\hspace{1cm} t \in \mathbb{N}_0.
\end{equation}
An important property of walks on {\it strongly connected} digraphs is (aperiodic) {\it ergodicity}.
A Markov chain (\ref{recursion}) is said to be {\it ergodic} if it exists 
$n_0 \in \{1,2,\ldots\} =\mathbb{N}$ such that
\begin{equation}
\label{ergodicity}
[\mathbf{P}(n)]_{ij} >0 , \hspace{1cm} \forall i,j =1,\dots,N \hspace{1cm} n\geq n_0, 
\end{equation}
is {\it strictly positive},
i.e.\ for each pair of nodes $i \to j$ there is at least one connecting path (and return path) not 
longer than $n_0$ jumps via directed connecting edges.
For a discussion of some aspects of ergodicity of biased walks in digraphs we refer to \cite{RiascosMichelitschPizarro2020} (and see also the references therein).
\\[1mm]
For later use we consider a Montroll-Weiss CTRW where a biased walk with transition matrix
(\ref{one-jump}) is subordinated to a counting process
with state-probabilities $\mathbb{P}(N(t)=n)=\Phi^{(n)}(t)$. Then the transition 
matrix (\ref{recursion}) with the elements $P_{ij}(t)$, indicating the probability to find 
the walker at time $t$ on node $j$ (with the given initial condition $P_{ij}(0)=\delta_{ij}$), is generalized by the Cox series \cite{Cox1967}
\begin{equation}
\label{coxseries}
[\mathbf{P}(t)]_{ij} = \sum_{n=0}^{\infty} \Phi^{(n)}(t) [\mathbf{W}^n]_{ij}, \hspace{1cm} t \in \mathbb{R}^{+}.
\end{equation}
We notice that the discrete-time Markov chain 
(\ref{recursion}) indeed is a special case of (\ref{coxseries}) when we account for a homogeneous
event stream with state probabilities $\Phi^{(n)}(t)= \Theta(t-n)-\Theta(t-(n+1))$, where $\Theta(t)$ is the Heaviside step function which is such that $\Theta(t)=1$  for $t\geq 0$ and 
$\Theta(t)=0$  for $t<0$ (in particular $\Theta(0)=1$). Indeed (\ref{coxseries}) boils down to $[\mathbf{P}(t)]_{ij}(t) = [\mathbf{W}^n]_{ij}$ and recovers then 
(\ref{recursion}) as $\Phi^{(n)}(t)=1$ for $t\in [n,n+1)$ and $\Phi^{(n)}(t)=0$ otherwise.

\section{BIASED WALKS WITH LONG-RANGE JUMPS}
\label{biasedwalks}
One goal of this paper is to introduce new types of biased walks with special emphasis on strictly increasing walks with long-range jumps
together with a systematic method for their construction. To this end, we consider a strongly connected (ergodic) digraph
with a Laplacian matrix $\mathbf{L}$ defined in (\ref{Laplaciandigraph}) characterizing its topology. 
Now we seek ``good'' Laplacian matrix functions $\mathbf{L} \rightarrow g(\mathbf{L})$
defining new topologies
such that the following ``good Laplacian properties'` are retained
\cite{RiascosMichelitsch-et-al2018,TMM-APR-ISTE2019,RiascosMichelitschPizarro2020}:
\\[3mm]
\noindent {\bf (i)} $\sum_{j=1}^N L_{ij}=0$ (i.e. the Laplacian matrix has a unique eigenvalue $\mu_1=0$ and a right eigenvector with constant (real) components 
$\langle i|v_1\rangle =1/\sqrt{N}$). This condition ensures that in a strongly connected graph the transition matrix has the unique real eigenvalue $\lambda_1=1$ of largest absolute value (reflecting the existence of a unique stationary distribution).
\\[3mm]
\noindent {\bf (ii)} $L_{ii} =k_i^{(\mathrm{out})} > 0$ ($\forall i=1,\ldots N$) tells us that each node has neighboring nodes, i.e.\ there are no isolated nodes. 
\\[3mm]
\noindent {\bf (iii)} $L_{ij} = -\Omega_{ij} \leq 0$ for $i\neq j$, i.e.\ the off-diagonal elements are all non-positive.
\\[2mm]
Clearly the conditions {\bf (i)-(iii)} maintain stochasticity of the 
transition matrix (\ref{one-jump}) and hence of (\ref{coxseries}).
\\[1mm]
We notice that the auxiliary Laplacian (\ref{reduced}) fulfills conditions {\bf (i)-(iii)} with uniform normalized degree ${\cal L}_{ii}=1$ $\forall i=1,\ldots N$. Without loss of generality it is now more convenient to work with the auxiliary Laplacian (\ref{reduced}) instead of the Laplacian matrix (\ref{Laplaciandigraph}).
We will prove subsequently that new non-trivial ``good'' Laplacian matrices retaining {\bf (i)-(iii)} are obtained by the class of matrix functions \cite{RiascosMichelitsch-et-al2018,TMM-APR-ISTE2019,RiascosMichelitschPizarro2020}
\begin{equation}
\label{bernsteinfunction}
g(\mathbf{\cal L}) = 
\int_0^{\infty} \left(\mathbf{1}-e^{-\tau\mathbf{\cal L}}\right)
\nu({\rm d}\tau)
\end{equation}
containing as argument the auxiliary Laplacian 
matrix $\mathbf{\cal L}$ (\ref{reduced}) and where ${\mathbf 1}$ denotes the unity matrix.
In the above integral, $\nu({\rm d}\tau)$ indicates a {\it L\'evy measure} for which (see e.g.\ \cite{Kochubei2011})
\begin{equation}
\label{Levyden}
\int_0^{\infty} \min\{1,\tau\}\nu({\rm d}\tau) < \infty 
\end{equation}
ensuring convergence of (\ref{bernsteinfunction}). 
Let us consider for a moment a scalar version $g(\mu)$ ($\mu \geq 0$) of (\ref{bernsteinfunction}). Then, we can define
a special class of Bernstein functions by the L\'evy-Khintchine representation
\begin{equation}
\label{bernstein}
g(\mu)= \int_0^{\infty} \left(1-e^{-\tau \mu}\right)
\nu({\rm d}\tau),
\end{equation}
which fulfills $g(\mu)|_{\mu=0}=0$. This property preserves the zero eigenvalue
in the matrix function (\ref{bernsteinfunction}). 
In subsequent applications we consider 
L\'evy measures $\nu({\rm d}\tau)= \nu(\tau){\rm d}\tau$ with non-negative
{\it L\'evy densities} $\nu(\tau)$.
The Laplacian Bernstein function (\ref{bernstein}) is also called Laplace exponent 
\cite{PolitoScalas2016}\footnote{This name comes from the fact that it occurs as a limit in strictly increasing subordinators $S_t$ as the Laplace transform $e^{-t g(\mu)} = \mathbb{E} e^{-\mu S_t}$ \cite{PolitoScalas2016}.}.
Generally, a non-negative function $f(a,b,\mu) \geq 0$ ($a,b \geq 0$) defined on $\mu \in [0,\infty)$ is said to be a {\it Bernstein function} if it fulfills $(-1)^{n-1}\frac{d^n}{d\mu^n}f(a,b,\mu) \geq 0$ ($n\geq 1$) and hence has the representation \cite{Kochubei2011}
\begin{equation}
\label{bernssteingen}
f_B(a,b,\mu)= a + b\mu + \int_0^{\infty} \left(1-e^{-\tau \mu}\right)
\nu({\rm d}\tau) .
\end{equation}
We consider Laplacian functions $g(\mu)$ (\ref{bernstein}) referring to the class of Bernstein functions 
with $a=b=0$, i.e.\ $g(\mu)=f_B(0,0,\mu)$. The condition $a=0$ is necessary to have $g(\mu)|_{\mu=0}=0$ for the sake of condition {\bf (i)}. The combination $a=0$, $b > 0$ gives an admissible Laplacian function fulfilling {\bf (i)}-{\bf (iii)}. However,
it adds a `trivial' contribution of the original Laplacian matrix ${\cal L}$. Therefore, we consider here only non-trivial
matrix functions (\ref{bernstein}) with also $b=0$.
We notice that the derivative $\frac{d}{d\mu}g(\mu)$ is a
completely monotonic function, i.e.
\begin{equation}
\label{completemono}
(-1)^{n-1} \frac{d^{n}}{d\mu^n} g(\mu) \geq 0 ,\hspace{1cm} n \in \mathbb{N}
\end{equation}
with $g(\mu) > 0$ for $\mu>0$ and $g(0)=0$. For theorems and a profound analysis of Bernstein and completely
monotonic functions and related subjects, see \cite{Widder1941, schilling}.
We now prove that (\ref{bernsteinfunction}) retains the properties {\bf (i)-(iii)}. For a proof in undirected networks 
with symmetric Laplacian matrices, see \cite{RiascosMichelitsch-et-al2018,TMM-APR-ISTE2019} 
and for the fractional Laplacian matrix function 
on digraphs consult \cite{RiascosMichelitschPizarro2020}.
We introduce a constant $\Lambda > 1$ such that the matrix
$\Lambda {\mathbf 1}-\mathbf{\cal L}$ has uniquely {\it non-negative} matrix elements, namely 
$(\Lambda-1)\delta_{ij}+\frac{\Omega_{ij}}{k^{(\mathrm{out})}_i} \geq 0$.
Then, for {\it ergodic} (i.e. strongly connected) digraphs it exists a 
$n_0>0$ such that all elements of the integer matrix power
\begin{equation}
\label{ergodicityequivalent}
[(\Lambda {\mathbf 1}-\mathbf{\cal L})^n]_{ij} > 0 ,\hspace{0.5cm} n\geq n_0,\hspace{0.5cm} \forall i,j=1,\dots N
\end{equation}
are {\it strictly positive}. Condition (\ref{ergodicityequivalent}) can be used as an equivalent definition of ergodicity of a digraph. One can further infer that if such a finite $n_0$ exists, then ergodic digraphs must have a finite number $N$ of states\footnote{The same holds also for undirected graphs \cite{TMM-APR-ISTE2019}.}.
Then, clearly the matrix exponential
\begin{equation}
\label{sumexpo}
e^{\Lambda {\mathbf 1}-\mathbf{\cal L}} =\sum_{n=0}^{\infty} 
\frac{1}{n!} [\Lambda {\mathbf 1}-\mathbf{\cal L}]^n
\end{equation}
is strictly positive\footnote{We call a matrix `{\it (strictly) positive}' if it has solely positive entries.}, as there are infinitely many uniquely positive matrices $\frac{1}{n!}[\Lambda {\mathbf 1}-\mathbf{\cal L}]^n$ contained in this series (namely those for $n \geq n_0$).
Since ${\bf 1}$ is commuting with any matrix $\mathbf{\cal L}$ we have
$e^{\Lambda {\mathbf 1}- \mathbf{\cal L}} = e^{\Lambda {\mathbf 1}} \cdot e^{-\mathbf{\cal L}}$
thus $[e^{\Lambda {\mathbf 1}-\mathbf{\cal L}}]_{ij} = e^{\Lambda} [e^{-\mathbf{\cal L}}]_{ij} >0$.
Now, since $e^{\Lambda} >0$ the matrix exponential $e^{-\mathbf{\cal L}}$ for an ergodic digraph is indeed strictly positive:
\begin{equation}
\label{uniquepos}
[e^{-\mathbf{\cal L}}]_{ij} >0 ,\hspace{1cm} \forall i,j=1,\ldots N.
\end{equation}
On the other hand, because of condition {\bf (i)} we have row-stochasticity
\begin{equation}
\label{sumrulenormalize}
\sum_{j=1}^N [e^{-\mathbf{\cal L}}]_{ij} =1 ,\hspace{1cm} \forall i,j=1,\ldots N
\end{equation}
As a consequence of eigenvalue zero of ${\cal L}$ this relation follows from $\sum_{j=1}^N[{\cal L}^n]_{ij} = \delta_{n0}$.
Then since all matrix elements (\ref{uniquepos}) are positive and row-normalized we have that
\begin{equation}
\label{diagooanl}
0< [e^{-\mathbf{\cal L}}]_{ij} < 1 ,\hspace{1cm} \forall i,j=1,\ldots N.
\end{equation}
Hence, it follows that $0< 1-[e^{-\mathbf{\cal L}}]_{ii} < 1$ thus the Bernstein matrix function 
$\mathbf{1}-e^{-\tau \mathbf{\cal L}}$ ($\tau>0$) retains the conditions {\bf (i)-(iii)} of a good Laplacian matrix.
By similar considerations for non-ergodic digraphs one can show that
the matrix exponential contains zero valued blocks of non-diagonal elements (where zero entries indicate the pairs of nodes where no finite directed paths exist). However, 
conditions {\bf (i)-(iii)} in the Bernstein matrix function (\ref{bernsteinfunction}) still are retained \cite{RiascosMichelitschPizarro2020}
\\[1mm]
Conditions {\bf (i)}-{\bf (iii)} remain also retained 
by the integration in (\ref{bernsteinfunction}) with a (non-negative) L\'evy measure
$\nu({\rm d}\tau)$ if the integral converges. This is definitely the case by virtue of (\ref{Levyden}) and {\it if and only if the eigenvalues of the auxiliary Laplacian ${\cal L}$ 
have solely non-negative real parts $\Re\{\mu_m\} \geq 0$}. 
This indeed is true as we show by the following brief proof. 
\\[1mm]
Denoting the eigenvalues of the auxiliary Laplacian (\ref{reduced}) by $\mu_m$ with $\mu_1=0$
and the eigenvalues of the one-step transition matrix (\ref{one-jump}) by $\lambda_m$
we have
\begin{equation}
\label{eigrel}
 \mu_m=1-\lambda_m  , \hspace{1cm} m=1,2,\ldots, N
\end{equation}
where $\mu_1=1-\lambda_1=0$ (thus $\lambda_1=1$ which is the unique so called Perron-Frobenius eigenvalue \cite{Frobenius1912}). Then, we have that $\lim_{n\rightarrow \infty}{\mathbf W}^n ={\mathbf W}^{\infty}$ remaining row-stochastic (corresponding to $\lambda_1=1$) which contradicts the existence of an exploding eigenvalue $|\lambda|^n\to \infty$. We can hence infer that for the complete set of complex eigenvalues it holds
\begin{equation}
\label{eigvals}
|\lambda_m| \leq 1 , \hspace{1cm} \forall m
\end{equation}
and with $\Re\{\lambda_m\} \leq |\lambda_m| \leq 1$ we have that
\begin{equation}
\label{largerzero}
\Re\{\mu_m\}=1-\Re\{\lambda_m\} \geq 1-|\lambda_m| \geq 0,  \hspace{1cm}  \forall m .
\end{equation}
We notice that in this proof we neither needed ergodicity nor that the digraph is finite, thus it also 
includes strictly increasing walks on the infinite integer line (see (\ref{eigenvaluescirculant})).
\\[1mm]
This concludes our proof that the L\'evy-Khintchine representation (\ref{bernsteinfunction}) indeed also converges 
for digraphs and is useful to construct good Laplacian matrix functions to define new biased walks with the one-step transition matrix

\begin{equation}
\label{allnodes}
W^{(g)}_{ij} =\delta_{ij}-\frac{g_{ij}(\cal L)}{g_{ii}(\cal L)}.
\end{equation}
It follows from above considerations (see especially (\ref{uniquepos})-(\ref{diagooanl}) with (\ref{bernsteinfunction})) for
{\it ergodic digraphs} that all off-diagonal elements of (\ref{allnodes}) are strictly positive $W^{(g)}_{ij} >0$
($i\neq j$) allowing the walker to reach any node in a single jump introducing new fully connected topologies.
In non-ergodic digraphs some off-diagonal elements take values null due to the absence of finite directed paths \cite{RiascosMichelitschPizarro2020}.
We also notice that as ${\cal L}$ is not symmetric and in some cases non-diagonisable having Jordan canonic forms 
\cite{Benzi-et-al-2020} where our proof remains valid also for these cases.
\section{CIRCULANT TRANSITION MATRICES AND STRICTLY INCREASING WALKS}
\label{DTRW}
In this section we consider the class of biased walks on the integer line 
with IID strictly positive integer increments (`jumps') $Z_j >0 $. Such a walk is defined by
%a random
%variable $Y_{n\in \mathbb{N}_0} \in \mathbb{N}_0$, namely
%
%
\begin{equation}
\label{strictly}
Y_n= \sum_{j=1}^n Z_j ,\hspace{1cm} Z_j \in \mathbb{N} ,\hspace{1cm} Y_0=0.
\end{equation}
The random walk (\ref{strictly}) is a discrete counterpart to a strict subordinator \cite{PachonPolitoRicciuti2018}.
There is a path
$a \to b$ only if $b>a$ but then no return path $b \to a$. Therefore in contrast to the walks on strongly connected (finite) digraphs, the walk
(\ref{strictly}) is not ergodic.
As an example let us recall the Sibuya walk. The Sibuya distribution has generating function \cite{PachonPolitoRicciuti2018}
\begin{equation}
\label{Sibyuagenerful}
{\bar w}_{\alpha}(u) = \mathbb{E}u^Z = 1-(1-u)^{\alpha} ,\hspace{0.5cm} \alpha \in (0,1),\hspace{0.5cm} |u| \leq 1.
\end{equation}
In a Sibuya walk the jumps $Z_j \in \mathbb{N}$ in (\ref{strictly}) follow the Sibuya distribution (also known under the name {\it Sibuya($\alpha$)}) which is defined 
for $\alpha \in (0,1)$ as
\begin{equation}
\label{Siuyaalpha}
\begin{array}{clc}
 \ds
\mathbb{P}(Z_j=k)=w_{\alpha}(k)= \frac{1}{k!}\frac{d^k}{du^{k}} {\bar w}_{\alpha}(u)\bigg|_{u=0} & =  \ds (-1)^{k-1}\frac{\alpha(\alpha-1)..(\alpha-k+1)}{k!} & \\ \\ \hspace{0.5cm}  \ds & = \ds (-1)^{k-1}\,\binom{\alpha}{k}=\frac{\alpha}{k}\frac{\Gamma(k-\alpha)}{\Gamma(1-\alpha)\Gamma(k)}    ,\hspace{0.3cm} k \in \mathbb{N} &
\end{array}
\end{equation}
with $w_{\alpha}(k)|_{k=0}= {\bar w}_{\alpha}(u)|_{u=0}=0$ (crucial for the occurrence 
of strictly positive jumps $Z_j >0$) and $w_{\alpha}(k) >0$ for all $k \in \mathbb{N}$.
In the limit $\alpha \to 1-$ the `trivial' distribution $w_{1}(k)=\delta_{1,k}$ is recovered (not called Sibuya)
where the walker makes unit jumps $Z_j=1$ to its right-sided neighbor node almost surely.
%The Sibuya walk is a discrete approximation of the stable strict subordinator \cite{PachonPolitoRicciuti2018}.
%
%
For later use we also mention the important property that Sibuya($\alpha$) is fat-tailed,
namely for $k$ large by using $\frac{\Gamma(k-\alpha)}{\Gamma(k)} \sim k^{-\alpha}$
we have the asymptotic behavior
\begin{equation}
\label{sibuya-asymptotics}
w_{\alpha}(k) \sim \frac{\alpha}{\Gamma(1-\alpha)}k^{-\alpha-1} = -\frac{k^{-\alpha-1}}{\Gamma(-\alpha)} >0 ,\hspace{1cm} (k \rightarrow \infty) .
\end{equation}
For the following analysis it will be convenient to utilize the connection of generating functions and the shift operators.
For a detailed outline we 
refer to our recent article \cite{MichelitschPolitoRiascos2020} and see also Appendix \ref{Appendix}.
\\[1mm]
Let us introduce the shift operator ${\hat T}_{a}$ ($a\in \mathbb{R}$) which is such that 
${\hat T}_{a}f(x)=f(x+a)$ and consider its {\it circulant} matrix representation
\begin{equation}
\label{circulantshift}
 {\hat T}_m f(p)= f(p+m) =: \sum_{q=-\infty}^{\infty}f(q)[{\hat T}_m]_{qp}=
 \sum_{q=-\infty}^{\infty}f(q)\delta_{q,p+m} ,\hspace{1cm} p,q, m \in \mathbb{Z}
 \end{equation}
and hence $[{\hat T}_m]_{qp}=\delta_{q,p+m}$.
This circulant structure remains true for all matrices of operator functions of shift operators.
All matrices $M_{mn}$ (we also write $M_{m,n}$) we are dealing with in the context of increasing walks
(\ref{strictly}) are characterized by the properties
\begin{equation}
\label{ciculant-matrices}
\begin{array}{clc}
\ds M_{m,n} = M_{m+s,n+s} & =  \ds M_{0,n-m} 
,\hspace{1cm} n, m, s \in \mathbb{Z} & \\ \\
\ds M_{m,n} = M_{0,n-m} & =  \ds 0 , \hspace{2cm}  n-m \in \mathbb{Z}_{-} &
\end{array}
\end{equation}
%
%
%i.e. $M_{m,n}$ are {\it upper triangular circulant} matrices, where
%all elements below the main diagonal are strictly null.
We call a matrix ``{\it upper triangular circulant}'' if it
fulfills (\ref{ciculant-matrices}),
i.e.\ all elements below the main diagonal are strictly null.
The Sibuya transition matrix (\ref{transitionoperator}) is a proto-typical example for such an upper triangular circulant matrix. 
We introduce the generating function of an upper triangular circulant matrix (\ref{ciculant-matrices}) 
as follows
\begin{equation}
\label{generatingmat}
\begin{array}{clc}
\ds 
{\bar M}(u)& \ds =  \sum_{k=-\infty}^{\infty}M_{0,k}u^k  = \sum_{k=0}^{\infty} M_{0,k}u^k, &  \ds |u|\leq 1 \\ \\ \ds 
M_{0,n} & \ds = \frac{1}{n!}\frac{d^n}{du^n} {\hat M}(u)|_{u=0} ,&  n \in \mathbb{N}_0 \\ \\
M_{0,n} & \ds = 0, & n < 0
\end{array}
\end{equation}
thus only non-negative powers $u^k$ ($k\geq 0$) are contained in ${\hat M}(u)$.
Then we observe the important property which holds for {\it upper triangular circulant matrices}
\begin{equation}
\label{property}
\begin{array}{clc}
\ds
\frac{1}{n!}\frac{d^n}{u^n}\left\{{\bar M}^{1}(u) {\bar M}^{2}(u)\right\}_{u=0} 
&  =  \ds  \sum_{k=0}^n M^{(2)}_{0,n-k}M^{(1)}_{0,k}  = 
\sum_{k=0}^n M^{(1)}_{0,k}M^{(2)}_{k,n}  =  \sum_{k=-\infty}^{\infty} M^{(1)}_{0,k}M^{(2)}_{k,n} & \\ \\
&= \ds  [{\mathbf M}^{(1)}{\mathbf M}^{(2)}]_{0,n} =  [{\mathbf M}^{(2)}{\mathbf M}^{(1)}]_{0,n} &
\end{array}
\end{equation}
connecting matrix multiplication and discrete 
(commuting) convolutions (See also \cite{MichelitschPolitoRiascos2020}).
\\[1mm]
The operator that emerges when replacing in (\ref{Sibyuagenerful}) 
$u$ with ${\hat T}_{-1}$ is hence the upper triangular circulant {\it one-step transition matrix} of the Sibuya walk \cite{RiascosMichelitschPizarro2020}, namely
\begin{equation}
\label{transitionoperator}
\begin{array}{clc} \ds
[{\bar w}_{\alpha}({\hat T}_{-1})]_{mn}=  [{\bar w}_{\alpha}({\hat T}_{-1})]_{0,n-m}& =& \ds 
[1-(1-{\hat T}_{-1})^{\alpha}]_{m,n} = \delta_{mn}-[{\cal L}^{\alpha}]_{m,n} , \hspace{0.5cm} \alpha \in (0,1] \\ \\ \ds
 &\hspace{-5cm}= &\hspace{-5cm} \ds \left\{ \begin{array}{l} \ds \sum_{k=1}^{\infty} (-1)^{k-1}\binom{\alpha}{k} \delta_{m,n-k} =\ds  (-1)^{n-m-1}\binom{\alpha}{n-m} ,\hspace{1cm} n-m > 0  \\ \\
 \ds 0 ,\hspace{1cm} n-m \leq 0 \end{array}\right.
\end{array}
\end{equation}
where $m,n \in \mathbb{Z}$. We used the following properties of the (unitary) shift operator with
$({\hat T}_{-1})^k ={\hat T}_{-k}$ and (\ref{circulantshift}) thus
\begin{equation}
\label{shiftprob}
[{\hat T}_{-k}]_{m,n} = \delta_{m,n-k} = \delta_{k,n-m} ,\hspace{1cm} m,n,k \in \mathbb{Z}
\end{equation}
is for $k\geq 0$ {\it upper} triangular circulant with entries $1$ in the $k$th upper side-diagonal ($n-m=k$), and $0$ elsewhere.
Formula (\ref{transitionoperator}) contains the upper triangular circulant fractional Laplacian matrix $g({\cal L})= {\cal L}^{\alpha}= (1-{\hat T}_{-1})^{\alpha}$ which has the elements
\begin{equation}
\label{matrixlapla}
[{\cal L}^{\alpha}]_{mn} = [(1-{\hat T}_{-1})^{\alpha}]_{mn}   =\left\{ \begin{array}{l} \ds
(-1)^{n-m}\binom{\alpha}{n-m} ,\hspace{1cm} n-m \in \mathbb{N}_0  \\ \\
 \ds 0 ,\hspace{1cm} n-m < 0 \end{array}\right.
\hspace{1cm} \alpha \in (0,1].
\end{equation}
For $\alpha=1$ the Laplacian matrix\footnote{From now on we refer the auxiliary Laplacian ${\cal L}$ to as Laplacian matrix.}
\begin{equation}
\label{Laplacian}
[{\cal L}]_{mn} = [(1-{\hat T}_{-1})]_{m,n}= \delta_{mn}-\delta_{m,n-1} ,\hspace{1cm} m,n \in \mathbb{Z}
\end{equation}
is recovered. It follows from Eq. (\ref{generatingmat}) that the generating function of Laplacian 
(\ref{Laplacian}) yields ${\bar {\cal L}}(u)=1-u$. Thus the generating function of the 
fractional Laplacian matrix (\ref{matrixlapla}) takes the form
${\bar {\cal L}^{\alpha}}(u)=(1-u)^{\alpha}$. The generating function of fractional Sibuya transition matrix (\ref{transitionoperator}) 
per construction recovers the generating function (\ref{Sibyuagenerful}) of Sibuya($\alpha$).
\\[1mm]
Further instructive for subsequent use is to consider the matrix elements of the matrix exponential
of the Laplacian (\ref{Laplacian}) which we conveniently obtain from its generating function in the form
\begin{equation}
\label{matrixexp}
[e^{-\tau{\cal L}}]_{0,n} = 
\frac{1}{n!}\frac{d^n}{du^n} e^{-\tau(1-u)}|_{u=0} = e^{-\tau}\frac{\tau^n}{n!} ,\hspace{1cm} \tau >0 ,\hspace{0.5cm} n\in {\mathbb N}_0.
\end{equation}
This matrix exponential is upper triangular circulant with a {\it Poisson distribution} in the non-vanishing entries. We directly verify in this representation the general properties (\ref{uniquepos})-(\ref{diagooanl})
of Laplacian matrix exponentials in digraphs. 
This relation underlines the utmost importance of the Poisson distribution
in strictly increasing walks.
\\[1mm]
These observations suggest that generating functions of good Laplacian
matrix functions $g(1-{\hat T}_{-1})$ of strictly increasing walks on the integer line such as occurring on right-hand side of (\ref{Cauchyproblemgen})
are with (\ref{bernsteinfunction}) obtained as
\begin{equation}
\label{LKbernstein}
g(1-u)= \int_{0}^{\infty}(1-e^{-\tau (1-u)})\nu({\rm d}\tau) = \sum_{k=0}^{\infty} [g(1-{\hat T}_{-1})]_{0,k}u^k ,\qquad |u| \leq 1 
\end{equation}
where clearly $g(1-{\hat T}_{-1})$ has {\it upper triangular circulant} matrix representation.
For instance for the L\'evy measure $\nu_{\alpha}({\rm d}\tau) = \frac{\alpha}{\Gamma(1-\alpha)}\tau^{-1-\alpha}{\rm d}\tau$ integral (\ref{LKbernstein}) converges within $\alpha \in (0,1)$ and yields generating function $(1-u)^{\alpha}$ of the Sibuya fractional Laplacian matrix (\ref{matrixlapla}).
The generating function of the one-step transition matrix then becomes with (\ref{LKbernstein}) and (\ref{allnodes})
\begin{equation}
\label{generfu}
\begin{array}{clc}
\ds {\bar W}^{(g)}(u) & = \ds  1-\frac{1}{g(1)}g(1-u)  ,\hspace{1cm} |u|\leq 1 ,\hspace{1cm} g(1)=[g({\cal L})]_{0,0}  &\\ \\
 & =\ds  \frac{1}{g(1)} \int_0^{\infty} e^{-\tau} (e^{\tau u}-1) \, \nu({\rm d}\tau)
\end{array}
\end{equation}
with the `generalized degree' $g(1)=g(1-u)|_{u=0}=\int_0^{\infty}[1-e^{-\tau}]\nu({\rm d}\tau) >0$.
The elements of the transition matrix are then straight-forwardly obtained as
\begin{equation}
\label{mteltran}
\begin{array}{clc}
\ds W^{(g)}_{0,n} 
& = \ds  \frac{1}{g(1)}\int_0^{\infty}\frac{\tau^n}{n!} e^{-\tau}\nu({\rm d}\tau) \,\, >0, & n \in \mathbb{N} \\ \\
\ds W^{(g)}_{0,0} &= 0. &
\end{array}
\end{equation}
The matrix elements are traced back
to Poisson terms $\frac{\tau^n}{n!} e^{-\tau}$ weighted by the L\'evy measure $\nu({\rm d}\tau)$ resulting in strictly positive
$W^{(g)}_{0,n} >0$ for $n>0$ and $W^{(g)}_{0,n}=0$ else introducing new topologies with directed edges allowing 
jumps of any positive integer size $Z_j \in \{1,2,\ldots\}=\mathbb{N}$.
\\[1mm]
For the class of L\'evy measures $\nu({\rm d}\tau)=\nu(\tau){\rm d}\tau $ with densities $\nu(\tau)$ which fulfill 
$\int_0^{\infty} \nu(\tau){\rm d}\tau < \infty$, i.e. having existing Laplace transforms 
\begin{equation}
\label{laplacelevy}
{\tilde \nu}(s) = 
\int_{0}^{\infty} e^{-\tau s} \nu(\tau){\rm d}\tau  < \infty ,\hspace{1cm} \Re\{s\} \geq 0
\end{equation}
we have $g(1-u) = {\tilde \nu}(0) - {\tilde \nu}(1-u)$ thus the generating function (\ref{generfu}) of the transition matrix writes
\begin{equation}
\label{generfultransmat}
\ds {\bar W}^{(g)}(u)= 
\frac{{\tilde \nu}(1-u)-{\tilde \nu}(1)}{{\tilde \nu}(0)-{\tilde \nu}(1)}.
\end{equation}
Due to the restriction of convergence of (\ref{laplacelevy}) we notice that relation
(\ref{generfu}) is more general than (\ref{generfultransmat}). 
\\[2mm]
After these considerations let us verify the crucial property (\ref{largerzero}) when we account for the (in our convention left-) eigenvectors of the unitary shift operator
$v_n=C e^{i\varphi n}$ ($\varphi \in (-\pi,\pi]$) with ${\mathbf v}\cdot{\cal L}=\sum_{m=-\infty}^{\infty} v_m [{\cal L}]_{mn}= (1-e^{-i\varphi})v_n$.
Hence the continuous complex eigenvalue spectrum of the Laplacian matrix (\ref{Laplacian}) is
\begin{equation}
\label{eigenvaluescirculant}
\mu(\varphi) =1-e^{-i\varphi} ,\hspace{1cm} \varphi \in (-\pi,\pi]
\end{equation}
with $\mu(\varphi)|_{\varphi=0}= 0$ and $\Re\{\mu(\varphi)\} = 1-\cos{\varphi} >0$ for $\varphi \neq 0$ ($\varphi \in (-\pi,\pi]$) in accordance with (\ref{largerzero})\footnote{Here ${\mathbf W}={\mathbf 1}-{\cal L}= [{\hat T}_{-1}]$ is unitary with $\lambda(\varphi)=e^{-i\varphi}$ thus for all eigenvalues holds $|\lambda(\varphi)|=1$.}.
Hence (\ref{bernsteinfunction}) converges and has the canonical representation
\begin{equation}
\label{spectral}
[g({\cal L})]_{n,m} = \frac{1}{2\pi} \int_{-\pi}^{\pi}e^{i\varphi(m-n)} g(1-e^{-i\varphi}){\rm d}\varphi
\end{equation}
where the eigenvalues of the Laplacian matrix function 
write with (\ref{bernsteinfunction})
\begin{equation}
\label{eigvallaplmathfu}
g(1-e^{-i\varphi}) = \int_0^{\infty}[1-e^{-\tau(1-e^{-i\varphi})}]\nu({\rm d}\tau)
\end{equation}
which is convergent since $\Re\{\mu(\varphi)\} = 1-\cos{\varphi} \geq 0$ ($\varphi \in (-\pi,\pi]$). Expanding
$e^{\tau e^{-i\varphi}}= \sum_{\ell=0}^{\infty}\frac{\tau^{\ell}}{\ell!}e^{-i\ell\varphi}$ leads to
\begin{equation}
\label{wearriveat}
\frac{1}{2\pi}\int_{-\pi}^{\pi}e^{i\tau e^{-i\varphi}}  e^{i(m-n)\varphi}{\rm d}\varphi 
=\left\{\begin{array}{l} \ds  \frac{\tau^{m-n}}{(m-n)!} \,\,,\hspace{1cm} m\geq n \\  \\ \ds  0 \,\, ,\hspace{2.5cm} m<n. \end{array}\right.
\end{equation}
Thus we get for the family of ``good Laplacian matrix functions'' as generators for strictly increasing walks the upper triangular circulant structure
\begin{equation}
\label{equival}
[g({\cal L})]_{n,m} = \left\{\begin{array}{clc} \ds 
\int_0^{\infty}\left[\delta_{n,m}-e^{-\tau}\frac{\tau^{m-n}}{(m-n)!}\right]\nu({\rm d}\tau) , & m-n \in \mathbb{N}_0 &\\ \\
\ds  0 , &  m-n <0  & \end{array}\right. 
\end{equation}
which clearly fulfills the ``good Laplacian properties'' {\bf (i)-(iii)} and is consistent with relation (\ref{mteltran}).
The convergence of these integrals is ensured by the property (\ref{Levyden}) of the L\'evy measure.
\\[3mm]
{\bf Mittag-Leffler transition matrix}
\\[2mm]
For later use we derive now in the framework of this approach the one-step transition matrix of a strictly
increasing walk on the integer line with jumps following a discrete approximation of the Mittag-Leffler density. 
Such a distribution was analyzed in
\cite{PilaiJayakumar1995} and generalizations in recent articles \cite{PachonPolitoRicciuti2018,MichelitschPolitoRiascos2020}.
To this end we consider a L\'evy measure $\nu_{ML,\alpha}({\rm d}\tau) = 
\nu_{ML,\alpha}(\tau){\rm d}\tau$ with L\'evy
density $\nu_{ML,\alpha}(\tau)= \lambda\tau^{\alpha-1}E_{\alpha,\alpha}(-\lambda \tau^{\alpha})$
which itself is a Mittag-Leffler density and has Laplace transform $\frac{\lambda}{\lambda+s^{\alpha}}$ ($\alpha \in (0,1]$, $\lambda >0$).
Then we get from the L\'evy-Khintchine representation (\ref{LKbernstein})
the Laplacian generating function 
\begin{equation}
\label{generfulaplacmittagleff}
g_{ML,\alpha}(1-u,\lambda) = 1-\frac{\lambda}{\lambda+(1-u)^{\alpha}} = \frac{(1-u)^{\alpha}}{\lambda+(1-u)^{\alpha}} ,\hspace{1cm} \lambda >0, \hspace{0.5cm} |u|\leq 1.
\end{equation}
Thus with (\ref{generfu}) the generating 
function of the ``{\it Mittag-Leffler transition matrix}'' yields
\begin{equation}\label{fracber}
{\bar W}_{ML,\alpha}(u,\lambda)= 1-\frac{1}{g_{ML,\alpha}(1,\lambda)}g_{ML,\alpha}(1-u,\lambda)   =\frac{\lambda}{\lambda+(1-u)^{\alpha}}(1-(1-u)^{\alpha}) ,\hspace{0.5cm} \alpha \in (0,1] \hspace{0.5cm} |u|\leq 1.
\end{equation}
This is the generating function of a discrete approximation of the {\it Mittag-Leffler} density
 introduced recently \cite{PachonPolitoRicciuti2018} in the context of discrete-time renewal processes and named there ``Discrete Mittag-Leffler'' (${\it DML}_A$) distribution (of so called ``type A'', see Remark 9 in that paper).
The elements of the {\it Mittag-Leffler transition matrix} are derived explicitly in relation (\ref{fracbernoulli-transmat}).
\\[1mm]
It can be seen beforehand that (\ref{fracber}) generates a discrete approximation of the Mittag-Leffler density by the following consideration.
By introducing the scaling $\lambda(h)=h^{\alpha}\lambda_0$ ($h,\lambda_0>0$ where $\lambda_0$ is a new positive constant independent of $h$) and with $u=e^{-hs}$ shows that the generating function (\ref{fracber}) converges to the Laplace transform of a Mittag-Leffler density, namely
\begin{equation}
\label{limit}
{\tilde {\cal W}}(\lambda_0,s) = \lim_{h\rightarrow 0} {\bar W}_{ML,\alpha}(e^{-hs},\lambda_0h^{\alpha})=\frac{\lambda_0}{\lambda_0+s^{\alpha}}
\end{equation}
where with ${\tilde f}(s)$ we denote the (spatial) Laplace transform of $f(x)$. It is instructive to
recall this limit in terms of shift operator representation (See Appendix \ref{Appendix}
and \cite{MichelitschPolitoRiascos2020} for details)
\begin{equation}
\label{shift-Mittag-Leffler}
\begin{array}{l}
\ds
{\cal W}_{ML,\alpha}(x,\lambda_0) = \lim_{h\to 0} {\bar W}_{ML,\alpha}({\hat T}_{-h},\lambda_0h^{\alpha})\delta_h(x) =
\lim_{h \to 0} \frac{\lambda_0h^{\alpha}}{\lambda_0h^{\alpha}+(1-{\hat T}_{-h})^{\alpha}}\delta_h(x)\\ \\ \ds =  \lim_{h \to 0} 
\frac{\lambda_0}{\lambda_0+h^{-\alpha} (1-{\hat T}_{-h})^{\alpha}}\delta_h(x)=    \frac{\lambda_0}{\lambda_0+D_x^{\alpha}}\delta(x)     = \lim_{h\to 0 } \frac{1}{h}[{\bar W}_{ML,\alpha}(\lambda_0h^{\alpha})]_{0,\frac{x}{h}}
\end{array}
\end{equation} 
which has Laplace transform converging to (\ref{limit}) and contains the {\it discrete $\delta$-distribution} $\delta_h(x)$ defined in (\ref{discrete-time-delta}).
Relation (\ref{shift-Mittag-Leffler}) defines the ``well-scaled'' continuous-space 
limit density kernel (`transition density kernel')
and converges to the Mittag-Leffler density (\ref{conti-space-lim}).
In the second line of (\ref{shift-Mittag-Leffler}) appears the (Riemann-Liouville) fractional derivative operator $D_x^{\alpha}$ of order $\alpha$ as 
$\lim_{h\to 0} h^{-\alpha}(1-e^{-hD_x})^{\alpha}=D_x^{\alpha}$ (with ${\hat T}_{-h}=e^{-hD_x}$).
\\[1mm]
Now we derive the matrix elements of the Mittag-Leffler transition matrix in explicit form.
For our convenience and later use we introduce the Pochhammer-symbol 
\begin{equation}
\label{Pochhammersymbol}
(c)_m =  
\frac{\Gamma(c+m)}{\Gamma(c)} = \left\{\begin{array}{l} 1 ,\hspace{1cm} m=0, \\ \\ 
                                                            
     c(c+1)\ldots (c+m-1) ,\hspace{1cm} m=1,2,\ldots \end{array}\right. 
\end{equation}
where especially $(c)_0=1$ and $(0)_{m}=\delta_{m0}$. Then we have with (\ref{fracber}) the expansions
\begin{equation}
\label{expandgenfu}
\begin{array}{clc} \ds 
\frac{\lambda}{\lambda+(1-u)^{\alpha}} & \ds  
=\left\{\begin{array}{clc} \ds 
\sum_{m=0}^{\infty} (-1)^m\lambda^{m+1} (1-u)^{-(m+1)\alpha}, 
& \ds |1-u|^{-1} \lambda^{\frac{1}{\alpha}} < 1 & \\ \\ \ds 
\sum_{m=0}^{\infty} (-1)^m\lambda^{-m}(1-u)^{\alpha m} , & \ds |1-u|^{-1} \lambda^{\frac{1}{\alpha}} > 1 &\end{array}\right. \\ \\
\ds {\bar W}_{ML,\alpha}(u,\lambda) & \ds = \left\{\begin{array}{clc}     \ds  
\sum_{m=0}^{\infty} (-1)^m\lambda^{m+1}
\left[(1-u)^{-(m+1)\alpha}- (1-u)^{-m\alpha}\right] 
,& \ds |1-u|^{-1} \lambda^{\frac{1}{\alpha}} < 1 & \\ \\ \ds
\sum_{m=0}^{\infty} (-1)^m\lambda^{-m}\left[(1-u)^{\alpha m}- (1-u)^{\alpha (m+1)}\right],
& \ds |1-u|^{-1} \lambda^{\frac{1}{\alpha}} > 1 & \end{array}\right.
\end{array}
\end{equation}
In order to determine the elements of the {\it Mittag-Leffler transition matrix}\footnote{As here all matrices commute we adopt the notation $\frac{\mathbf A}{\mathbf B} = {\mathbf A}\cdot 
{\mathbf B}^{-1}$.} we have to account for the cases of convergence in (\ref{expandgenfu}) at $u=0$ to arrive at
\begin{equation}
\label{fracbernoulli-transmat}
\begin{array}{clc}
\ds [W_{ML,\alpha}(\lambda)]_{0,n} & \ds = 
\left[\frac{\lambda{\mathbf 1}}{\lambda{\mathbf 1}+{\cal L}^{\alpha}}({\mathbf 1}-{\cal L}^{\alpha})\right]_{0,n} =
\frac{1}{n!}\frac{d^n}{du^u}{\bar W}_{ML,\alpha}(u,\lambda)|_{u=0} , \hspace{3.5cm} n \in \mathbb{N}_0 &\\ \\ 
& \ds  = \left\{\begin{array}{clc} \ds \frac{\lambda}{n!} \sum_{m=0}^{\infty} (-1)^m\lambda^{m} \left[ \frac{\Gamma(\alpha (m+1)+n)}{\Gamma(\alpha(m+1))} - 
\frac{\Gamma(\alpha m +n)}{\Gamma(\alpha m)}\right] ,  &  \ds 0 < \lambda <1
 & \\ \\ \ds 
\frac{(-1)^n}{n!}\sum_{m=0}^{\infty} (-1)^m\lambda^{-m} \left[\frac{\Gamma(\alpha m+1)}{\Gamma(\alpha m-n+1)}- \frac{\Gamma(\alpha (m+1) +1)}{\Gamma(\alpha (m+1)-n+1)}\right] ,& \ds  \lambda >1 \end{array}\right.
 \end{array}
\end{equation}
where these series converge absolutely by accounting for the asymptotic behavior of the
terms for $m$ large as $\frac{\Gamma(\alpha (m+1)+n)}{\Gamma(\alpha(m+1))} \sim m^n \lambda^m$ for $\lambda <1$, and
in the same way $\sim m^n \lambda^{-m} $ for $\lambda>1$, respectively. We also have $[W_{ML,\alpha}(\lambda)]_{0,0}=0$.
\\[1mm]
It is worthy also to consider the case $\alpha=1$ where (\ref{fracber}) takes the form
\begin{equation}
\label{casealpha1}
\begin{array}{clc} 
\ds 
{\bar W}_{ML,1}(u,\lambda)& \ds  = \frac{ u \lambda}{\lambda +1 - u}  
  = \frac{\lambda}{(\lambda+1)} \frac{u}{(1-\frac{u}{\lambda+1})} &\\ \\ & = \ds \frac{\lambda}{(\lambda+1)} 
\sum_{k=0}^{\infty} \frac{u^{k+1}}{(\lambda+1)^k}  = \sum_{k=1}^{\infty} pq^{k-1} u^k  & \\ \\  & \ds p=\frac{\lambda}{\lambda+1} ,\hspace{1cm} q= \frac{1}{\lambda+1} , & p+q=1
\end{array}
\end{equation} 
and hence we obtain for $\alpha=1$ the transition matrix
\begin{equation}
\label{alspha1}
\begin{array}{clc}
 \ds [W_{ML,1}(\lambda)]_{0,n}& = pq^{n-1} ,\hspace{1cm} n \in \mathbb{N} & \\ \\ \ds
 [W_{ML,1}(\lambda)]_{0,0} & = 0 &
 \end{array}
 \end{equation}
which recovers the {\it geometrical distribution} (See also \cite{PachonPolitoRicciuti2018}).
By using the explicit formula (\ref{fracbernoulli-transmat}) for the Mittag-Leffler transition matrix, it is now not a big deal to perform explicitly the “well-scaled” continuous-space limit (\ref{shift-Mittag-Leffler}).
Accounting for the asymptotic relation $ \frac{(\beta)_n}{n!} \sim \frac{n^{\beta-1}}{\Gamma(\beta)}$ ($n=\frac{x}{h}\rightarrow \infty$), and by introducing the scaling 
$\lambda(h)=\lambda_0h^{\alpha} \rightarrow 0$ which is covered in relation (\ref{fracbernoulli-transmat}) by the case $0<\lambda<1$, we arrive at
\begin{equation}
\label{MLdensity}
\begin{array}{clc} \ds 
{\cal W}_{ML,\alpha}(x,\lambda_0)&= \ds  \lim_{h\to 0} \frac{1}{h}[{\bar W}_{ML,\alpha}(\lambda_0h^{\alpha})]_{0,\frac{x}{h}}, \hspace{1cm}
x \in h\mathbb{N}_0 &
\\ \\ & = \ds  \lim_{h\to 0}
\sum_{m=0}^{\infty} (-1)^m\lambda_0^{m+1} 
\frac{1}{\frac{x}{h}!} h^{\alpha(m+1)-1}
\left[((m+1)\alpha)_{\frac{x}{h}}-(m\alpha)_{\frac{x}{h}}\right]. &  
\end{array}
\end{equation}
We employ here the Pochhammer symbol (\ref{Pochhammersymbol}) and we account for
\begin{equation}
\label{limits}
\lim_{h\to 0} h^{\alpha(m+1)-1}((m+1)\alpha)_{\frac{x}{h}} = 
\lim_{h\to 0} h^{\alpha(m+1)-1} \frac{1}{\Gamma[\alpha(m+1)]}\left(\frac{x}{h}\right)^{\alpha(m+1)-1} =\frac{x^{\alpha (m+1)-1}}{\Gamma[\alpha(m+1)]}
\end{equation}
with $x \in \mathbb{R}^+$.
The second term in (\ref{MLdensity}) has a vanishing limit
\begin{equation}
\label{whereas}
\lim_{h\to  0} h^{\alpha(m+1)-1} (m\alpha)_{\frac{x}{h}} =
\lim_{h\to 0}  h^{\alpha(m+1)-1}
\frac{1}{\Gamma(\alpha m)} \left(\frac{x}{h}\right)^{\alpha m-1} = \lim_{h\to 0} 
h^{\alpha}\frac{x^{\alpha m-1}}{\Gamma(\alpha m)} = 0 .
\end{equation}
We hence obtain for the well-scaled continuous-space limit (\ref{MLdensity}) as anticipated the Mittag-Leffler density
\begin{equation}
\label{conti-space-lim}
\begin{array}{clc}
\ds {\cal W}_{ML,\alpha}(x,\lambda_0) &= \ds  \lambda_0x^{\alpha-1} 
\sum_{m=0}^{\infty} \frac{(-\lambda_0 x^{\alpha})^m}{\Gamma(\alpha m+\alpha)} & \\ \\  &\ds  = 
\lambda_0 x^{\alpha-1}E_{\alpha,\alpha}(-\lambda_0 x^{\alpha}) 
=\frac{d}{dx}[1-E_{\alpha}(-\lambda_0x^{\alpha})] , & \hspace{0.5cm} x\in \mathbb{R}^{+} ,\hspace{0.5cm} \alpha \in (0,1]
\end{array}
\end{equation}
containing the 
generalized Mittag-Leffler function 
$E_{\alpha,\gamma}(z)=\sum_{m=0}^{\infty}\frac{z^m}{\Gamma(\alpha m+\gamma)}$. 
For $\alpha=1$ (\ref{conti-space-lim}) recovers the exponential density 
${\cal W}_{ML,1}(x,\lambda_0) = \lambda_0 e^{-\lambda_0 x}$ which also is the continuous-space limit 
of (\ref{alspha1}) obtained with $p(h)= \frac{\lambda_0h}{\lambda_0 h+1}$ and
$q(h)= \frac{1}{\lambda_0h+1}$ thus
\begin{equation}
\label{contilimalpha1}
\begin{array}{clc}
\ds {\cal W}_{ML,1}(x,\lambda_0) &= \ds 
\lim_{h\rightarrow 0} \frac{1}{h}[{\bar W}_{ML,1}(\lambda_0h)]_{0,\frac{x}{h}} & \\ \\ 
& \ds = \lim_{h\rightarrow 0}
\lambda_0  (1+ \lambda_0h)^{-\frac{x}{h}} & \\ \\ & =  \ds \lambda_0 e^{-\lambda_0 x} = \frac{d}{dx}(1-e^{-\lambda_0 x}), &\hspace{1cm}  x\in \mathbb{R}^{+}.
\end{array}
\end{equation}
\section{SPACE-TIME GENERALIZATIONS OF THE POISSON PROCESS}
\label{space-time-gens}
\subsection{The classical cases}
During the last two decades, an increasing interest in generalizations of the Poisson renewal process has emerged. The most natural generalization probably is obtained when the exponential waiting time density is
generalized by a Mittag-Leffler density. The resulting generalization is the
fractional Poisson process which was introduced and analyzed by several authors
\cite{RepinSaichev2000,Laskin2003,MainardiGorenfloScalas2004,BeghinOrsinger2009,GorenfloMainardi2013}.
Then space-time generalizations of the Poisson process were introduced
such as the `{\it space-fractional Poisson Process}', the `{\it space-time fractional Poisson process}' \cite{OrshingerPolito2012} and further generalizations of the latter 
\cite{PolitoScalas2016} were developed within the last decade. Generally these space-time generalizations can be seen as strictly increasing walks on the integer line time-changed with an independent renewal process.
Before we introduce in Section \ref{space-time-frac} such a generalization, let us briefly recall these
classical cases.
\subsection{Poisson process}
We consider the Cauchy problem
\begin{equation}
\label{cauchypoisson}
\begin{array}{clc} 
 \ds \frac{d}{dt}p_n(t) &= \ds -\xi p_n(t)+\xi p_{n-1}(t)= -\xi (1-{\hat T}_{-1})p_n(t)  
 ,& \ds p_n(0) = \ds\delta_{n0} ,\hspace{1cm}  \xi >0 ,\, t\geq 0 \\ \\ \ds
 \frac{d}{dt} \mathbf{p}(t) & =  \ds -\xi \mathbf {p}(t)\cdot \mathbf{\cal L}. &
 \end{array}
\end{equation}
Indeed Cauchy problem (\ref{cauchypoisson}) defines the time-evolution of the state probabilities in the standard Poisson counting process
which is the simplest variant of the general class (\ref{Cauchyproblemgen}) (where ${\cal L}$ is the circulant Laplacian matrix
(\ref{Laplacian})). Equation (\ref{cauchypoisson}) is the Kolmogorov-Feller (also referred to as Kolmogorov-forward) equation which is solved by the state probabilities of the standard Poisson counting process $p_n(t) = \mathbb{P}({\cal N}(t)=n)$ ($n\in \mathbb{N}_0$, $t\geq 0$) where ${\cal N}(t) \in \mathbb{N}_0$ counts the events within the time interval $[0,t]$.
We can conceive the Poisson counting process in the Montroll-Weiss CTRW picture \cite{MontrollWeiss1965} as a strictly increasing random walk 
with one-step transition matrix 
$[{\mathbf W}]_{m,n}=[{\hat T}_{-1}]_{m,n}=\delta_{m,n-1}$ subordinated to a Poisson process where at each arrival the walker makes a jump of increment $Z_j=1$ almost surely.
Clearly the solution of the Cauchy problem (\ref{cauchypoisson}) is the well-known standard Poisson distribution
\begin{equation}
\label{standarxPoisson-dis}
p_n(t) =  [\mathbf{p}(0) \cdot e^{-\xi t\mathbf{\cal L}}]_{n} =
\frac{1}{n!}\frac{d^n}{du^n}e^{-\xi t(1-u)}|_{u=0}= \frac{(\xi t)^n}{n!}e^{-\xi t} , \hspace{0.5cm} n \in \mathbb{N}_0 , \hspace{0.5cm} t \in \mathbb{R}^{+}.
\end{equation}
\subsection{Fractional Poisson process}
The most natural time-generalization of the Poisson process is defined by the Cauchy problem
\begin{equation}
\label{difffract} 
\begin{array}{clc} \ds \frac{d^{\beta}}{dt^{\beta}}p_n(t) & =  \ds -\xi (1-{\hat T}_{-1})p_n(t) , \hspace{1cm} \xi>0 ,\hspace{1cm} p_n(0)=\delta_{n0}, \hspace{0.5cm} \beta \in (0,1] & \\ \\ \ds 
\frac{d^{\beta}}{dt^{\beta}} \mathbf{p}(t) & = \ds -\xi \mathbf {p}(t)\cdot \mathbf{\cal L} &
\end{array}
\end{equation}
which is the well-known fractional Kolmogorov-Feller equation of the {\it fractional Poisson process} \cite{Laskin2003} where
$\frac{d^{\beta}}{dt^{\beta}}$ denotes the {\it Caputo fractional derivative} of order $\beta$ defined as, e.g. 
\cite{SamkoKilbasMarichev1993}
\begin{equation}
\label{caputo}
\frac{d^{\beta}}{d t^{\beta}}p(t)= \frac{1}{\Gamma(1-\beta)}\int_0^t(t-\tau)^{-\beta}\frac{d}{d\tau}p(\tau){\rm d}\tau 
,\hspace{1cm} \beta \in (0,1]
\end{equation}
and in the limit 
$\lim_{\beta\rightarrow 1-} \frac{(t-\tau)^{-\beta}}{\Gamma(1-\beta)}=\delta(t-\tau)$ 
(\ref{caputo}) recovers the first-order derivative $\frac{d}{dt}p(t)$ 
thus the fractional Poisson process turns into the standard Poisson process.
The solution of (\ref{difffract}) writes
\begin{equation}
\label{Cauchysolutiofrac}
\mathbf{p}(t)= \mathbf{p}(0)\cdot E_{\beta}(-\xi t^{\beta}\mathbf{\cal L}) 
= [\ldots , [E_{\beta}(-\xi t^{\beta}\mathbf{\cal L})]_{0,n}, \ldots]
\end{equation}
where comes into play the Mittag-Leffler matrix function
\begin{equation}
\label{Mittag-Leffler-matfu}
E_{\beta}(-\xi t^{\beta}\mathbf{\cal L}) = 
\sum_{m=0}^\infty\frac{(-\xi t^{\beta})^m}{\Gamma(\beta m+1)} \mathbf{\cal L}^m 
\end{equation}
with the Mittag-Leffler function $E_{\beta}(z) = \sum_{m=0}^\infty\frac{z^m}{\Gamma(\beta m+1)}$.
Taking into account the generating function $\sum_{k=0}^{\infty}[\mathbf{\cal L}^m]_{0,k}u^k=(1-u)^m$, 
we obtain for the generating function of the Mittag-Leffler matrix function (\ref{Mittag-Leffler-matfu}) the expression $\sum_{k=0}^{\infty}[E_{\beta}(-\xi t^{\beta}\mathbf{\cal L})]_{0,k}u^k =E_{\beta}(-\xi t^{\beta}(1-u))$. Hence with (\ref{generatingmat}) we get for the components of the state vector (\ref{Cauchysolutiofrac})
\begin{equation}
\label{weget}
\begin{array}{clc}\ds 
\ds p_n^{\beta}(t) &= \ds  [E_{\beta}(-\xi t^{\beta}\mathbf{\cal L})]_{0,n} =
\frac{1}{n!}\frac{d^n}{du^n}E_{\beta}(- (1-u)\xi t^{\beta})|_{u=0}, \hspace{0.5cm} n \in \mathbb{N}_0 , \hspace{0.5cm} t \in \mathbb{R}^{+} &
\\ \\ \ds & = \ds 
\frac{(\xi t^{\beta})^n}{n!}\sum_{m=0}^{\infty}\frac{(m+n)!}{m!}
\frac{(-\xi t^{\beta})^m}{\Gamma((m+n)\beta +1)}.&
\end{array}   
\end{equation}
We identify (\ref{weget}) indeed with the state probabilities of Laskin's fractional Poisson distribution \cite{Laskin2003} which also was
derived in different manners by several further authors \cite{RepinSaichev2000,MainardiGorenfloScalas2004,BeghinOrsinger2009}.
For $\beta=1$ (\ref{weget}) recovers the Poisson distribution (\ref{standarxPoisson-dis}).
\subsection{ Space-time fractional Poisson process}
Orsingher and Polito have generalized the Cauchy problems (\ref{cauchypoisson}), (\ref{difffract}) by introducing a spatial generalization $\mathbf{\cal L} \rightarrow \mathbf{\cal L}^{\alpha}$ where the right-hand sides of (\ref{cauchypoisson}) and (\ref{difffract}), respectively, are generalized by the (Sibuya-) fractional Laplacian matrix (\ref{matrixlapla}). They named these processes 
{\it space-fractional Poisson process} and {\it space-time fractional Poisson process}, respectively \cite{OrshingerPolito2012}. 
Indeed the space-fractional Poisson process is a Sibuya-walk subordinated to an independent Poisson process, and the
space-time fractional Poisson process is a Sibuya walk time-changed with an independent fractional Poisson process.
The Cauchy problem defining the space-time fractional Poisson process then writes
\begin{equation}
\label{space-time-fract} 
\begin{array}{clc} \ds \frac{d^{\beta}}{dt^{\beta}}p_n(t) & = \ds -\xi (1-{\hat T}_{-1})^{\alpha} p_n(t) ,\hspace{1cm} \alpha,\beta \in (0,1] ,\hspace{0.2cm} \xi>0,  \hspace{0.2cm} t\geq 0 & \\ \\
\ds p_n(t)|_{t=0} &= \ds  \delta_{n0}  & \\ \\ \ds 
\frac{d^{\beta}}{dt^{\beta}} \mathbf{p}(t)& = -\xi \mathbf {p}(t)\cdot \mathbf{\cal L}^{\alpha} &
\end{array} 
\end{equation}
For $\beta=1$ the space-fractional Poisson process is recovered and for $\alpha=1$ the fractional Poisson process.
We obtain the state vector solving (\ref{space-time-fract}) as
\begin{equation}
\label{state-vector-space-time-frac}
{\mathbf{p}}^{\alpha,\beta}(t) = {\mathbf{p}}(0)\cdot E_{\beta}(-\xi t^{\beta}{\cal L}^{\alpha}) ,\hspace{1cm} \alpha,\beta \in (0,1] ,\hspace{0.2cm} t\geq 0 
\end{equation}
having generating function
\begin{equation}
\label{generfuspace-time-frac}
{\bar p}^{\alpha,\beta}(u,t) =  E_{\beta}(-\xi t^{\beta}(1-u)^{\alpha})
\end{equation}
where ${\bar p}^{\alpha,\beta}(u,t)|_{u=1}=1$ reflects row-stochasticity of the upper triangular 
circulant transition matrix $E_{\beta}(-\xi t^{\beta}{\cal L}^{\alpha})$ (normalization of the state probabilities). 
The state distribution is obtained as
\begin{equation}
\label{thespace-time-frac}
\begin{array}{l}
 \ds p_n^{\alpha,\beta}(t) = \frac{1}{n!}\frac{d^n}{du^n}E_{\beta}(- (1-u)^{\alpha}\xi t^{\beta})|_{u=0} ,\hspace{0.5cm}  n \in \mathbb{N}_0 ,\hspace{0.5cm} t\geq 0\\ \\ \ds
\hspace{0.5cm} = \sum_{m=0}^{\infty}\frac{(-\xi t^{\beta})^m}{\Gamma(\beta m+1)} \frac{1}{n!} \frac{d^n}{du^n} (1-u)^{\alpha m}|_{u=0} = \frac{(-1)^n}{n!} \sum_{m=0}^{\infty} \frac{\Gamma(\alpha m+1)}{\Gamma(\alpha m +1-n)} \frac{(-\xi t^{\beta})^m}{\Gamma(\beta m+1)}
\end{array}
\end{equation}
which is the expression obtained by Orsingher and Polito (Eq. (2.28) in \cite{OrshingerPolito2012})
and recovers for $\beta=1$ their expression for the {\it space-fractional Poisson process} derived in the same paper \cite{OrshingerPolito2012}. Further generalizations of the
space-fractional Poisson are considered in the references \cite{toa,gar}.
\subsection{Well-scaled diffusion limits}
We consider briefly the diffusion limit of standard Poisson $\alpha,\beta=1$. To this end we define a well-scaled continuous-space 
limit in (\ref{standarxPoisson-dis}) by introducing the scaling $\xi(h)=\xi_0h^{-1}$ (where $\xi_0>0$ is a new dimensional constant independent of $h$) to arrive at
\begin{equation}
\label{specialcase11}
{\cal P}_{1,1}(x,t) = \lim_{h\rightarrow 0} e^{-\xi_0 t \, h^{-1}(1-e^{-hD_x})}\delta_h(x)
= e^{-\xi_0 t D_x}\delta(x)
={\hat T}_{-\xi_0t}\,\delta(x) = \delta(x-\xi_0 t)
\end{equation}
which is a moving Dirac $\delta$-distribution propagating with constant velocity $\xi_0$ in the positive $x$-direction. It is immediately checked that (\ref{specialcase11}) solves the continuous-space limit of (\ref{cauchypoisson}) which is the Cauchy problem 
\begin{equation}
\label{Cspacelim}
\frac{\partial}{\partial t}{\cal P}_{1,1}(x,t) =-\xi_0 \frac{\partial}{\partial x} {\cal P}_{1,1}(x,t), \qquad
{\cal P}_{1,1}(x,t)|_{t=0}=\delta(x) .
\end{equation}
\subsection{ Diffusion limit of space-time fractional Poisson}
Further consider the space-time fractional Poisson process $\alpha \in (0,1)$ and $\beta \in (0,1]$ with the scaling $\xi(h)=\xi_0h^{-\alpha}$ ($\xi_0>0$).
Then we can write the continuous-space diffusion equation
which emerges from the well-scaled limit of (\ref{space-time-fract}).
By accounting for the limit $\lim_{h\rightarrow 0}h^{-\alpha}(1-{\hat T}_{-h})^{\alpha} =D_x^{\alpha}$
which takes the Riemann-Liouville fractional derivative of order $\alpha$ (e.g. \cite{SamkoKilbasMarichev1993,Podlubny1999,michelCFM2011} and many others) we obtain
for the well-scaled limit of (\ref{space-time-fract}) the space-time fractional diffusion equation
\begin{equation}
\label{Cspacelimalpha}
\frac{\partial^{\beta}}{\partial t^{\beta}}{\cal P}(x,t) =-\xi_0 D_x^{\alpha}{\cal P}(x,t), \qquad 
{\cal P}(x,t)|_{t=0}=\delta(x).
\end{equation}
On the right hand side of (\ref{Cspacelimalpha}) occurs the spatial Riemann-Liouville fractional derivative of order $\alpha \in (0,1)$ defined as, e.g. \cite{SamkoKilbasMarichev1993,Podlubny1999}
\begin{equation}
\label{RLalpha}
D_x^{\alpha}{\cal P}(x,t) = \frac{\partial}{\partial x}\int_0^x
\frac{(x-\tau)^{-\alpha}}{\Gamma(1-\alpha)} {\cal P}(\tau,t){\rm d}\tau ,\hspace{1cm} \alpha \in (0,1).
\end{equation}
We notice that the fractional Laplacian 
${\cal L}^{\alpha}=[(1-{\hat T}_{-1})^{\alpha}]$ on the right-hand side of (\ref{space-time-fract})
does not contain an own scaling parameter. Therefore, in order to get an existing limit we have to rescale the constant $\xi(h)=\xi_0h^{-\alpha}$ (having dimension $[\sec]^{-\beta}$) where $\xi_0$ has units $[\sec^{-\beta}\text{cm}^{\alpha}]$.
\section{SPACE-TIME MITTAG-LEFFLER PROCESS}
\label{space-time-frac}
Having recalled these classical cases, 
we introduce here a generalization of the Poisson process which is a strictly increasing walk on the integer line with {\it Mittag-Leffler jumps}
taking place at independent fractional Poisson arrival times. We also highlight the connections with the Montroll-Weiss CTRW picture in more details.
We call this process
`{\it space-time Mittag-Leffler process}' which we define by a Cauchy problem of the general type (\ref{Cauchyproblemgen}), namely
\begin{equation}
\label{Space-time-fractional-Bernoulli}
\begin{array}{lcl}
 \ds
\frac{d^{\beta}}{dt^{\beta}}p_n(t) &= &\ds   - \xi(\lambda+1)\frac{(1-{\hat T}_{-1})^{\alpha}}{(\lambda+(1-{\hat T}_{-1})^{\alpha})}p_n(t)  ,\hspace{1cm} \xi,\, \lambda \, >\, 0, \hspace{0.5cm} t \geq 0 ,\hspace{0.5cm} \alpha,\, \beta \in (0,1] , \hspace{0.5cm} n \in \mathbb{N}_0
\\ \\  \ds p_n(t)|_{t=0} & = & \ds \delta_{0,n}    \\ \\ \ds
\frac{d^{\beta}}{dt^{\beta}}{\mathbf p}(t)& = & \ds -\frac{\xi}{[g_{ML,\alpha}]_{0,0}} {\mathbf p}(t)\cdot g_{ML,\alpha}({\cal L}) , \hspace{0.5cm} [{\mathbf p}(t)\cdot g_{ML,\alpha}({\cal L})]_n= \sum_{m=0}^n p_m(t) [g_{ML,\alpha}({\cal L})]_{0,n-m}.
\end{array} 
\end{equation}
The right-hand side contains the good Laplacian matrix function $g_{ML,\alpha}({\cal L}) =\frac{{\cal L}^{\alpha}}{\lambda+{\cal L}^{\alpha}}$ of (\ref{generfulaplacmittagleff}) (with generalized degree $[g_{ML,\alpha}]_{0,0}=g_{ML,\alpha}(1-u,\lambda)|_{u=0}=\frac{1}{\lambda+1}$) and
generates discrete Mittag-Leffler jumps with transition matrix (\ref{fracbernoulli-transmat}). 
The last line indicates the matrix representation and uses the upper triangular circulant property of $g_{ML,\alpha}({\cal L})$
(see (\ref{ciculant-matrices})) and
$\frac{d^{\beta}}{dt^{\beta}}$ stands for the {\it Caputo fractional derivative} (\ref{caputo}).
The generating function representation of Cauchy problem
(\ref{Space-time-fractional-Bernoulli}) then writes 
\begin{equation}
\label{writethegenfunction}
\begin{array}{lcl}
 \ds
\frac{d^{\beta}}{dt^{\beta}}{\bar p}(u,t) & = 
& \ds -\xi \frac{(\lambda+1)(1-u)^{\alpha}}{\lambda +(1-u)^{\alpha}}{\bar p}(u,t) ,\hspace{1cm} |u| \leq 1, \\ \\ \ds \ds
{\bar p}(u,t)|_{t=0} & = &1.
\end{array}
\end{equation}
For the solution of (\ref{Space-time-fractional-Bernoulli}) we can write
\begin{equation}
\label{soltionspace-timefrac}
{\mathbf p}^{\alpha,\beta}_{\lambda,\xi}(t) = {\mathbf p}(0)\cdot 
E_{\beta}\left(-\xi t^{\beta} (\lambda+1)\frac{{\cal L}^{\alpha}}{\lambda{\mathbf 1}+{\cal L}^{\alpha}}\right)
\end{equation}
with the generating function of the state-probabilities
\begin{equation}
\label{generfufracber}
\begin{array}{lcl}
 \ds
{\bar p}^{\alpha,\beta}_{\lambda,\xi}(u,t) & =  & \ds
E_{\beta}\left(\frac{-\xi t^{\beta} (\lambda+1)(1-u)^{\alpha} }{\lambda+(1-u)^{\alpha}}\right) \\ \\ &=& \ds \sum_{m=0}^{\infty}
 \frac{[-\xi t^{\beta} (\lambda+1)]^m}{\Gamma(\beta m+1)}\frac{(1-u)^{\alpha m}}{[\lambda+(1-u)^{\alpha}]^{m}}
\end{array}
 \end{equation}
involving the standard Mittag-Leffler function $E_{\beta}(z)$. In the Poisson limit $\beta=1$ due to
$E_1(z)= e^{z}$ the Mittag-Leffler functions in all relations recover exponentials and 
the Caputo derivative recovers the standard first order time-derivative.
We observe that ${\bar p}^{\alpha,\beta}_{\lambda,\xi}(u,t)|_{u=1} = 1$ reflecting normalization of the state probabilities.
Expression (\ref{generfufracber}) contains 
\begin{equation}
\label{discrtePrabhakar}
\begin{array}{lcl}
\ds
{\bar {\cal E}}_{\alpha}^{(m)}(\lambda,u) & = & \ds [g_{ML,\alpha}(1-u)]^m = \frac{1}{[1+\lambda(1-u)^{-\alpha}]^m} \\ \\ 
\ds & =  & \ds  \sum_{s=0}^{\infty} \frac{(m)_s}{s!}(-\lambda)^s(1-u)^{-\alpha s},\hspace{1cm} \lambda^{\frac{1}{\alpha}}|1-u|^{-1} < 1
\end{array}
\end{equation}
and plainly
\begin{equation}
\label{othercase}
\begin{array}{lcl}
\ds {\bar {\cal E}}_{\alpha}^{(m)}(\lambda,u) & = & \ds  [g_{ML,\alpha}(1-u)]^m = \frac{\lambda^{-m}(1-u)^{\alpha m}}{(1+\lambda^{-1}(1-u)^{\alpha})^m} \\ \\ & = & 
\ds
 \sum_{s=0}^{\infty} \frac{(m)_s}{s!}(-1)^s\lambda^{-s-m}(1-u)^{\alpha (s+m)} 
,\hspace{1cm} \lambda^{\frac{1}{\alpha}}|1-u|^{-1} > 1
\end{array}
\end{equation}
where $(m)_s=m(m+1)..(m+s-1)$ stands for the Pochhammer-symbol (\ref{Pochhammersymbol}).
${\bar {\cal E}}_{\alpha}^{(m)}(\lambda,u)$
is the generating function of a discrete version of a so called Prabhakar kernel \cite{Prabhakar1971,Giusti2020,MichelitschPolitoRiascos2020}. 
We will show this connection explicitly in subsequent analysis of well-scaled continuous-space limits.
The state-probabilities solving (\ref{Space-time-fractional-Bernoulli})
are then with (\ref{generfufracber}) 
and (\ref{discrtePrabhakar}), (\ref{othercase}) obtained as
\begin{equation}
\label{state-vector}
 \begin{array}{lcl}
 \ds
 p^{\alpha,\beta}_{\lambda,\xi,n}(t)& =& \ds \frac{1}{n!}\frac{d^n}{du^n} 
 E_{\beta}\left(\frac{- (\lambda+1)(1-u)^{\alpha} \xi t^{\beta}}{\lambda+(1-u)^{\alpha}}\right)|_{u=0} ,\hspace{1cm} n \in \mathbb{N}_0  \\ \\ & = & \ds
 \sum_{m=0}^{\infty}
 \frac{[- (\lambda+1)\xi t^{\beta}]^m}{\Gamma(\beta m+1)} {\cal E}_{\alpha}^{(m)}(\lambda,n).
 \end{array}
\end{equation}
Since we take the derivatives at $u=0$ we need to account for that (\ref{discrtePrabhakar})
converges at $u=0$ for $0< \lambda <1$, whereas (\ref{othercase}) converges at $u=0$ for 
$\lambda >1$. We hence arrive at
\begin{equation}\label{prabhakarapprox}
 {\cal E}_{\alpha}^{(m)}(\lambda,n) = \left\{\begin{array}{clc} \ds 
 \frac{1}{n!} \sum_{s=0}^{\infty} \frac{(-\lambda)^s (m)_s}{s!} \frac{\Gamma(\alpha s+n)}{\Gamma(\alpha s)} , & 0<\lambda <1\\ \\ 
 \ds   \frac{(-1)^n\lambda^{-m}}{n!} \sum_{s=0}^{\infty} \frac{(-1)^s \lambda^{-s} (m)_s}{s!} \frac{\Gamma(\alpha(s + m)+1)}{\Gamma(\alpha (s+m)-n+1)} ,&  \lambda > 1
 \end{array} \right.
 \hspace{1cm} n \in \mathbb{N}_0.
\end{equation}
One can easily verify in view of the asymptotic scaling $\frac{\Gamma(s\alpha+a)}{\Gamma(\alpha s +b)} \sim s^{a-b}$ and $\frac{(m)_s}{s!}\sim s^{m-1}$ for $s \to \infty$ that the series (\ref{prabhakarapprox})
for the two cases converge absolutely.
\begin{figure}[!t] 
\begin{center}
\includegraphics*[width=0.9\textwidth]{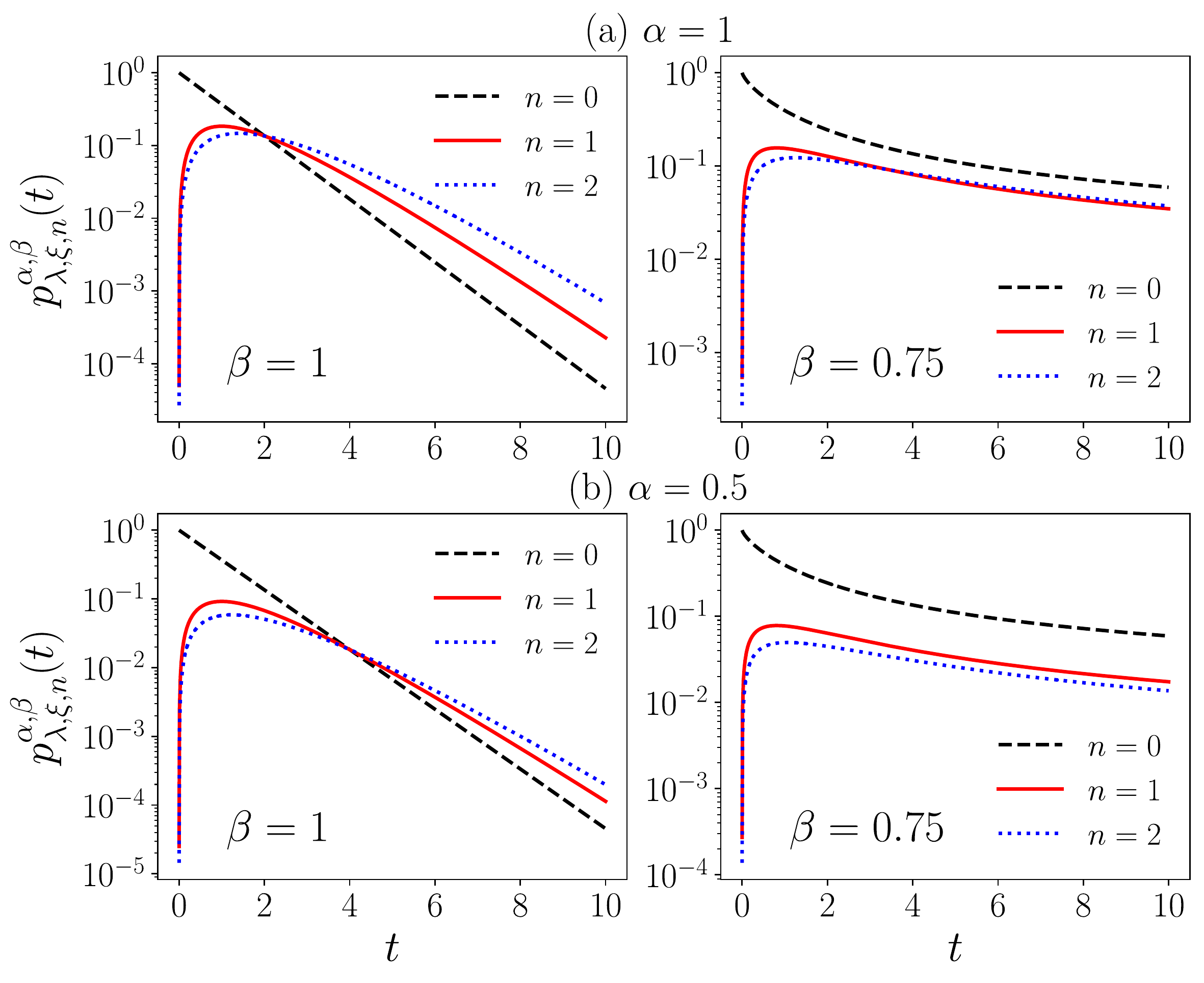}
\end{center}
\vspace{-7mm}
\caption{\label{Fig1} State-probabilities $p^{\alpha,\beta}_{\lambda,\xi,n}(t)$ of the space-time Mittag-Leffler process as a function of $t$ for: (a) $\alpha=1.0$ and (b) $\alpha=0.5$.  The results are obtained numerically using Eq. (\ref{state-vector}) for the values $n=0,1,2$ and $\beta=1$ (in the left panels) and $\beta=0.75$ (presented in the right panels); we maintain constant the parameters $\lambda=1$, $\xi=1$. 
} 
\end{figure}
We directly verify the initial condition $p^{\alpha,\beta}_{\lambda,\xi,n}(t)|_{t=0}={\cal E}_{\alpha}^{(0)}(\lambda,n) = \delta_{0,n}$.
For $n=0$ we have ${\cal E}_{\alpha}^{(m)}(\lambda,0) =(1+\lambda)^{-m}$ thus
the `survival probability' (probability that the walker at time $t$ still is in the initial state $n=0$)
is Mittag-Leffler, namely 
\begin{equation}
\label{survival}
 p^{\alpha,\beta}_{\lambda,\xi,0}(t) = {\bar p}^{\alpha,\beta}_{\lambda,\xi}(u,t)|_{u=0}    = E_{\beta}(-\xi t^{\beta})
 \end{equation}
which is necessarily the survival probability in the fractional Poisson process. We notice that for $\beta \in (0,1)$ the space-time Mittag-Leffler process is non-Markovian with long memory features, whereas for $\beta=1$ it becomes Markovian
due to the memoryless nature of the standard Poisson process, e.g.
\cite{Laskin2003,SaichevZaslavsky1997}. 
\\[1mm]
In Figure \ref{Fig1} the time-dependence of the state probabilities (\ref{state-vector}) for $n=0,1,2$ is depicted for different values of $\alpha,\beta \in (0,1]$, respectively. 
The survival probability $n=0$ is independent of $\alpha$ and given by the Mittag-Leffler survival probability (\ref{survival}) which turns for
$\beta=1$ into an exponential (see the plots on the left) with initial conditions
$p^{\alpha,\beta}_{\lambda,\xi,0}(t)|_{t=0}=1$ whereas for $n=1,2$ the zero initial conditions can be seen in the plots.
For large dimensionless times $\tau = t\xi^{\frac{1}{\beta}} \to \infty$ we have for $\beta \in (0,1)$ a power-law decay
$p^{\alpha,\beta}_{\lambda,\xi,n}(t) \sim \tau^{-\beta}$ (see asymptotic relation (\ref{state-prob-as})).
\\[2mm]
\noindent {\bf Asymptotic behavior}
\\[2mm]
It is worthwhile to consider the asymptotic behavior for $n$ large and finite (dimensionless) times
$\xi^{\frac{1}{\beta}}t$. To this end consider generating function for $u\to 1-0$ for $\alpha \in (0,1)$, namely
\begin{equation}
\label{expandML} 
 E_{\beta}\left(\frac{-(\lambda+1)\xi t^{\beta}(1-u)^{\alpha}}{\lambda(1+\lambda^{-1}(1-u)^{\alpha})}\right)
= 1 -\frac{(\lambda+1) \xi t^{\beta}}{\Gamma(\beta+1)\lambda}(1-u)^{\alpha} +{\cal O}(1-u)^{\alpha}
\end{equation}
where ${\cal O} (1-u)^{\alpha}$ indicates higher orders such that
$\lim_{u\to 1-}(1-u)^{-\alpha} {\cal O} (1-u)^{\alpha}=0$. We notice in view of (\ref{expandML})
that this asymptotic behavior is for $\alpha \in (0,1)$ of the same fat-tailed type as for the state probabilities (\ref{thespace-time-frac}) of the space-time fractional Poisson process and also as Sibuya($\alpha$), namely
\begin{equation}
\label{largenasym}
\begin{array}{l}
\ds
p^{\alpha,\beta}_{\lambda,\xi,n}(t) \sim -\frac{(\lambda+1)\xi t^{\beta}}{\lambda \Gamma(\beta+1)}
\frac{1}{n!}\frac{d^n}{du^n}(1-u)^{\alpha}|_{u=0} \\ \\ \ds \hspace{0.5cm} =\frac{(\lambda+1)\xi t^{\beta}}{\Gamma(\beta+1)}
(-1)^{n-1}\binom{a}{n} 
\sim \frac{\alpha(\lambda+1) \xi t^{\beta}  n^{-\alpha-1}}{\lambda \Gamma(1-\alpha)\Gamma(1+\beta)} ,\hspace{1cm} n\rightarrow \infty
\end{array} 
\end{equation}
where $\alpha \in (0,1)$ and $\beta \in (0,1]$. The fat-tailed behavior $p^{\alpha,\beta}_{\lambda,\xi,n}(t) \sim n^{-\alpha-1}$ reflects the occurrence of long-range forward jumps and is equivalent to the divergence of the first moment, namely
$\frac{d}{du} {\bar p}^{\alpha,\beta}_{\lambda,\xi,n}(u,t)|_{u=1}= \sum_{n=0}^{\infty} p^{\alpha,\beta}_ {\lambda,\xi,n}(t)  n \rightarrow \infty$ for $\alpha\in (0,1)$.
The asymptotic form (\ref{largenasym}) contains also the `well-scaled' continuous-space limit ($h\to 0$: $\lambda(h)=\lambda_0h^{\alpha}$ and $x \in h\mathbb{N}_0 \to \mathbb{R}^{+}$), namely 
\begin{equation}
\label{continlim}
{\cal P}^{\alpha,\beta}_{\lambda_0,\xi}(x,t) = \lim_{h\to 0} \frac{1}{h} p^{\alpha,\beta}_{\lambda_0h^{\alpha},\xi,\frac{x}{h}}(t) \sim \frac{\alpha x^{-\alpha-1} \, \xi t^{\beta} }{\lambda_0 \Gamma(1-\alpha)\Gamma(1+\beta)}, \hspace{0.5cm} t\geq 0 ,\hspace{0.5cm} x\rightarrow \infty .
\end{equation}
The ``well-scaled'' continuous-space limiting procedure will be justified in subsequent paragraph in more details. $\lambda_0>0$ is a new constant (independent of $h$) and has units $[\text{cm}]^{-\alpha}$. 
The constant $\xi>0$ has physical dimension $[\sec]^{-\beta}$ thus (\ref{continlim}) is a spatial density of units $[\text{cm}]^{-1}$.
\\[1mm]
Then let us also consider the asymptotic behavior for large (dimensionless) time $t\xi^{\frac{1}{\beta}} \to \infty$ and finite $n$.
From the asymptotic behavior of the scalar Mittag-Leffler function 
$E_{\beta}(-at^{\beta}) \sim \frac{1}{a}\frac{t^{-\beta}}{\Gamma(1-\beta)}$ with 
$a=\frac{\xi(\lambda+1)}{1+\lambda(1-u)^{-\alpha}}$ follows for the asymptotic behavior of the generating function
(\ref{generfufracber})
\begin{equation}
\label{generas}
{\bar p}^{\alpha,\beta}_{\lambda,\xi}(u,t) \sim \frac{(1+\lambda(1-u)^{-\alpha})}{(\lambda+1)\xi}\frac{t^{-\beta}}{\Gamma(1-\beta)} ,
\hspace{1cm} t\xi^{\frac{1}{\beta}} \to \infty ,\hspace{1cm} \beta \in (0,1), \hspace{0.5cm} \alpha \in (0,1]
\end{equation}
and hence for the state probabilities
\begin{equation}
\label{state-prob-as}
p^{\alpha,\beta}_{\lambda,\xi,n}(t) \sim  
\frac{(\delta_{n,0}+\lambda\frac{(\alpha)_n}{n!})}{(\lambda+1)}\frac{t^{-\beta}}{\xi \Gamma(1-\beta)},
\hspace{1cm} t\xi^{\frac{1}{\beta}} \to \infty , n\in \mathbb{N}_0 ,\hspace{0.5cm} \beta \in (0,1), \hspace{0.5cm} \alpha \in (0,1]
\end{equation}
where for $n=0$ the pure Mittag-Leffler asymptotic behavior
of the survival probability (\ref{survival}) $p^{\alpha,\beta}_{\lambda,\xi,0}(t)\sim \frac{t^{-\beta}}{\xi \Gamma(1-\beta)}$
is obtained. We also get the asymptotic behavior in view of (\ref{generas}) for the well-scaled continuous-space limit by
($h\to 0$: $\lambda(h)=\lambda_0h^{\alpha}$ and $x \in h\mathbb{N}_0 \to \mathbb{R}^{+}$) for large dimensionless times and finite continuous state variable $x \in [0,\infty)$
\begin{equation}
\label{contilim}
\begin{array}{clc}
\ds {\cal P}^{\alpha,\beta}_{\lambda_0,\xi}(x,t) & \ds  \sim 
 \lim_{h\to 0} \frac{1}{h} p^{\alpha,\beta}_{\lambda_0h^{\alpha},\xi,\frac{x}{h}}(t) 
= \lim_{h \to 0} \frac{\left(\delta_h(x)+\lambda_0h^{\alpha-1}\frac{(\alpha)_{\frac{x}{h}}}{\frac{x}{h}!}\right)}{(\lambda_0h^{\alpha}+1)}
\frac{t^{-\beta}}{\xi \Gamma(1-\beta)} & \\ \\
 \ds {\cal P}^{\alpha,\beta}_{\lambda_0,\xi}(x,t) & \ds \sim \left(\delta(x) + \lambda_0\frac{x^{\alpha-1}}{\Gamma(\alpha)}\right) \frac{t^{-\beta}}{\xi\Gamma(1-\beta)} , \hspace{1cm} t\xi^{\frac{1}{\beta}} \to \infty, \hspace{0.5cm} x \in \mathbb{R}^{+}, \hspace{0.5cm}\beta \in (0,1) &
\end{array}
\end{equation}
where the asymptotic $t^{-\beta}$-power-law decay reflects the non-Markovian long memory feature 
of the process.
\\[2mm]
\noindent {\bf Well-scaled continuous-space limit}
\\[2mm]
Having derived the state probabilities (\ref{state-vector}) solving Cauchy problem
(\ref{Space-time-fractional-Bernoulli}), we are now interested in the well-scaled continuous-space limit density solving a (forward) diffusion equation which turns out to refer to the general class (\ref{general}).
We can define this diffusion limit by the scaling assumption $\lambda(h)=h^{\alpha}\lambda_0$ where $\lambda_0>0$ is independent of $h$. Here the constant $\xi$ does not need to be rescaled in order to obtain an existing limit.
We define the `well-scaled' continuous-space limit state density kernel by 
\begin{equation}
\label{well-scaled-limit}
\begin{array}{l}
 \ds 
{\cal P}_{\lambda_0,\xi}^{\alpha,\beta}(x,t) = \lim_{h\rightarrow 0}E_{\beta}\left(\frac{-\xi(\lambda_0h^{\alpha}+1) t^{\beta}}{1+ \lambda_0h^{\alpha}(1-{\hat T}_{-h})^{-\alpha}}\right)\delta_h(x)
=\lim_{h\rightarrow 0} \frac{1}{h}
p^{\alpha,\beta}_{\lambda_0h^{\alpha},\xi,\frac{x}{h}}(t) \\ \\
\ds \hspace{0.5cm} = E_{\beta}\left(\frac{-\xi t^{\beta}}{1+\lambda_0D_x^{-\alpha}}\right)\delta(x)
\end{array}
\end{equation}
where we employed the discrete-$\delta$ distribution $\delta_h(x)$ defined in (\ref{discrete-time-delta}) and its limiting behavior
(\ref{dirac-delta}) and in this limiting process $x\in h\mathbb{N}_0 \to \mathbb{R}^{+}$.
$D_x^{-\alpha}$ indicates the Riemann-Liouville fractional integral operator
of order $\alpha$. By taking into account (\ref{state-vector}) we get
\begin{equation}
\label{explicitexpspacelimit}
{\cal P}_{\lambda_0,\xi}^{\alpha,\beta}(x,t) = \lim_{h\rightarrow 0} 
 \sum_{m=0}^{\infty}
 \frac{[- (\lambda_0h^{\alpha}+1)\xi t^{\beta}]^m}{\Gamma(\beta m+1)} \frac{1}{h} {\cal E}_{\alpha}^{(m)}\left(\lambda_0h^{\alpha},\frac{x}{h}\right)
\end{equation}
where we use the asymptotic relation $\frac{(\beta)_n}{n!} \sim 
\frac{n^{\beta-1}}{\Gamma(\beta)}$. We have to consider in (\ref{prabhakarapprox}) the case $0<\lambda(h)<1$ as $\lambda(h)=\lambda_0h^{\alpha} \to 0$ to evaluate the continuous-space limit
\begin{equation}
\label{toevaluate}
\begin{array}{clc}
\ds e_{\alpha,0}^{m}(-\lambda_0,x) &= \ds 
    \lim_{h\to 0} \frac{1}{h} {\cal E}_{\alpha}^{(m)}\Bigl(\lambda_0h^{\alpha}, \frac{x}{h}\Bigr) &\\ \\  &=  \ds \lim_{h\rightarrow 0} \left[1+\lambda_0h^{\alpha}(1-{\hat T}_{-h})^{-\alpha}\right]^{-m}\delta_h(x) & \\ \\ 
     &=
 \ds \lim_{h\to 0} \sum_{s=0}^{\infty} \frac{(m)_s(-\lambda_0)^s}{s!}h^{\alpha s} \sum_{n=0}^{\infty}
 \frac{(\alpha s)_n}{n!}\delta_h(x-hn)  & \\ \\ 
& =  \ds \lim_{h\rightarrow 0} \delta_h(x) + \sum_{s=1}^{\infty} \frac{(m)_s(-\lambda_0)^s}{s!}h^{\alpha s-1} 
 \frac{(\alpha s)_{\frac{x}{h}}}{\frac{x}{h}!} & \\ \\
 \ds   e_{\alpha,0}^{m}(-\lambda_0,x) & =  \ds \delta(x) + 
 \sum_{s=1}^{\infty}  \frac{ (m)_s (-\lambda_0)^sx^{\alpha s-1} }{s!\Gamma(\alpha s)} , & m=0,1,2,.. \in \mathbb{N}_0 , \hspace{0.5cm} x \in  \mathbb{R}^{+}
 \\ \\ \ds  & = \ds D_x\left[\Theta(x)E_{\alpha,1}^m(-\lambda_0x)\right]  & 
\end{array}
\end{equation}
and for $m=0$ we have $e_{\alpha,0}^0(-\lambda_0,x)=\lim_{h\rightarrow 0}\delta_h(x)=\delta(x)$. In the last line we have
used the property of the Heaviside step function $D_x\Theta(x)=\delta(x)$ and it comes into play the so called
{\it Prabhakar function} which is a generalization of the Mittag-Leffler function introduced by
Prabhakar \cite{Prabhakar1971} as
\begin{equation}
\label{prabhakarfu}
E_{\alpha,\beta}^{\gamma}(z)= \sum_{s=0}^{\infty}\frac{(\gamma)_s}{s!}\frac{z^s}{\Gamma(\alpha s+\beta)} ,\hspace{0.5cm}
\alpha,\beta,\gamma \in {\mathbb C},\hspace{0.5cm} \Re\{\alpha\} > 0.
\end{equation}
We used in the deduction (\ref{toevaluate}) that ${\hat T}_{-hn}\delta_h(x)= \delta_h(x-hn)=\frac{1}{h}\delta_{\frac{x}{h},n}$
where $h\rightarrow 0: \,  x \in h \mathbb{N}_0 \rightarrow \mathbb{R}^{+}$
(see also Appendix \ref{Appendix}, Eq. (\ref{discrete-time-delta})).
The continuous-space limit expression (\ref{toevaluate}) can be identified with a so called {\it Prabhakar kernel} $ e_{\alpha,0}^m(-\lambda_0,x)$. The Prabhakar kernel $e_{\alpha,\beta}^{\gamma}(\zeta,x)$ was introduced
by Giusti \cite{Giusti2020} (where we use here his definition) as the kernel with Laplace transform ${\tilde e}_{\alpha,\beta}^{\gamma}(\zeta,x)(s)= s^{-\beta}(1-\zeta s^{-\alpha})^{-\gamma}$.
It follows that (\ref{prabhakarapprox}) indeed is a discrete approximation of Prabhakar kernel (\ref{toevaluate}).
With (\ref{toevaluate}) and (\ref{explicitexpspacelimit}) we obtain for 
the {\it state density kernel}
\begin{figure}[!t] 
\begin{center}
\includegraphics*[width=0.8\textwidth]{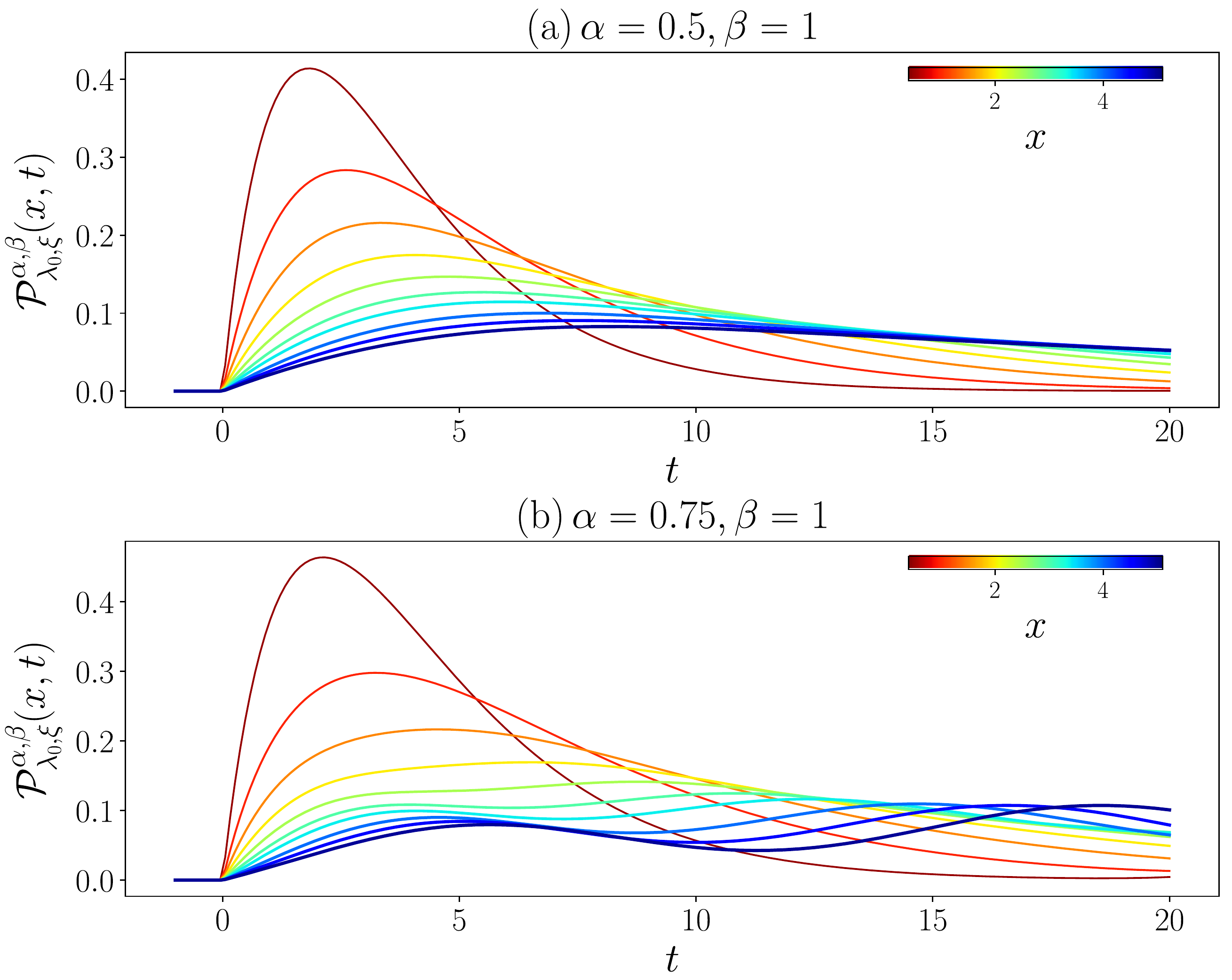}
\end{center}
\vspace{-4mm}
\caption{\label{Fig2}  State density kernel $\mathcal{P}_{\lambda_0,\xi}^{\alpha,\beta}(x,t)$ as a function of $t$. The results are obtained numerically from Eq. (\ref{finalresultspacecontilimit}) for: (a) $\alpha=0.5$ and (b) $\alpha=0.75$ maintaining constant $\lambda_0=1$, $\xi=1$ and $\beta=1$. We present the values for $x=0.5,1,1.5,2,\ldots, 5$ with different colors codified in the colorbar.
} 
\end{figure}
\begin{equation}
\label{finalresultspacecontilimit}
\begin{array}{clc} \ds 
{\cal P}_{\lambda_0,\xi}^{\alpha,\beta}(x,t) & = \ds 
\sum_{m=0}^{\infty}
 \frac{(-\xi t^{\beta})^m}{\Gamma(\beta m+1)}  e_{\alpha,0}^m(-\lambda_0,x), \hspace{1cm} \lambda_0,\xi >0 ,\hspace{0.5cm} \alpha,\beta \in (0,1], \hspace{0.5cm} x,t \geq 0 & \\ \\
 &= \ds  D_x\left[\Theta(x) \Theta(t) \sum_{m=0}^{\infty}
 \frac{(-\xi t^{\beta})^m}{\Gamma(\beta m+1)} E_{\alpha,1}^m(-\lambda_0x) \right] &
 \end{array}
\end{equation}
which has units $[cm]^{-1}$. We directly verify the initial
condition ${\cal P}_{\lambda_0,\xi}^{\alpha,\beta}(x,t)|_{t=0}= e_{\alpha,0}^0(-\lambda_0,x) =\delta(x)$.
\\[1mm]
The time-dependence of the state density kernel (\ref{finalresultspacecontilimit}) is plotted in Figure \ref{Fig2}
for the Poisson limit $\beta=1$ for two different values of $\alpha$.
Increasing values of the state variable $x$ are indicated by colors turning from red (small $x>0$) to blue (large $x$).
We observe that for larger $\alpha$ in the lower plot
the state density exhibits increasingly oscillating behavior for increasing values of $x$. We will come back to this issue subsequently
when we discuss the emerging continuous-space limit diffusion equation governing the time-evolution of 
the state density.
\\[1mm]
It is now only a small step to derive the continuous-space forward diffusion equation
of generalized fractional type which is solved by the state density kernel (\ref{finalresultspacecontilimit}). To this end
we deduce the convolution kernel of continuous-space limit of the right-hand side
in (\ref{Space-time-fractional-Bernoulli}), by the well-scaled limit (see also relations (\ref{generfulaplacmittagleff}) and (\ref{generfu}))
\begin{equation}
\begin{array}{clc}
\ds 
{\cal G}_{\xi,\lambda_0,\alpha}(x) &= \ds \lim_{h\rightarrow 0} \frac{g_{ML,\alpha}(1-{\hat T}_{-h},\lambda_0h^{\alpha})}{g_{ML,\alpha}(1,\lambda_0h^{\alpha})}\delta_h(x)   =     \lim_{h\rightarrow 0} 
\frac{(\lambda_0h^{\alpha}+1)}{1+\lambda_0h^{\alpha}(1-{\hat T}_{-h})^{-\alpha}}\delta_h(x) & \\ & =  \ds 
\frac{1}{(1+\lambda_0D_x^{-\alpha})}\delta(x) = \sum_{m=0}^{\infty} (-\lambda_0)^mD_x^{-m\alpha}\delta(x).
\end{array}
\end{equation}
Therefore
\begin{equation}
\label{space-contilim}
\begin{array}{clc}
\ds 
{\cal G}_{\xi,\lambda_0,\alpha}(x) &=  \ds 
D_x \frac{D_x^{\alpha-1}}{D_x^{\alpha}+\lambda_0}\delta(x) = \ds D_x[\Theta(x)E_{\alpha}(-\lambda_0x^{\alpha})]  &
\\ \\ 
 &= \ds  \delta(x)+ \Theta(x) D_x E_{\alpha}(-\lambda_0x^{\alpha}) = \delta(x) -\Theta(x)\lambda_0x^{\alpha-1}E_{\alpha,\alpha}(-\lambda_0x^{\alpha}) & \\ \\
 & = \ds \delta(x) - {\cal W}_{ML,\alpha}(x,\lambda_0). 
\end{array}
\end{equation}
We call this kernel `{\it Laplacian density}'.
This result also is obtained by taking into account its spatial Laplace transform 
$\frac{s^{\alpha}}{s^{\alpha}+\lambda_0}$ (obtained from the scaling limit of the first line in (\ref{space-contilim})).
We notice that the Laplacian density still maintains the `distributional versions' of the good Laplacian properties {\bf (i)-(iii)}: We have
$\int_0^{\infty}{\cal G}_{\xi,\lambda_0,\alpha}(x){\rm d}x=0$ corresponding to {\bf (i)}, ${\cal G}_{\xi,\lambda_0,\alpha}(x)<0 $ for $x>0$
(condition {\bf (iii)}), and $\lim_{\epsilon\to 0+}\int_0^{\epsilon}{\cal G}_{\xi,\lambda_0,\alpha}(x){\rm d}x=1$ corresponds to {\bf (ii)}. 
In the last line of (\ref{space-contilim}) we account for the Mittag-Leffler transition density kernel ${\cal W}_{ML,\alpha}(x,\lambda_0)=
\lambda_0x^{\alpha-1}E_{\alpha,\alpha}(-\lambda_0x^{\alpha})$ (\ref{conti-space-lim}) occurring as continuous-space limit 
of the Mittag-Leffler matrix (\ref{fracbernoulli-transmat}).
The continuous-space-time Cauchy problem then writes
\begin{equation}
\label{Cauchycont}
\begin{array}{clc}
\ds
\frac{\partial^{\beta}}{\partial t^{\beta}}{\cal P}(x,t) &= \ds -\xi \int_{-\infty}^{\infty}{\cal P}(\tau,t)
{\cal G}_{\xi,\lambda_0,\alpha}(x-\tau){\rm d}\tau ,
\hspace{0.5cm} x,t \geq 0 ,  \hspace{0.5cm} \alpha, \beta \in (0,1]  &
\\ \\  & = \ds -\xi \int_{0}^{x}{\cal P}(\tau,t)
{\cal G}_{\xi,\lambda_0,\alpha}(x-\tau){\rm d}\tau & \\ \\
& = \ds -\xi {\cal P}(x,t) + \xi \int_0^x {\cal W}_{ML,\alpha}(x-\tau,\lambda_0){\cal P}(\tau,t){\rm d}\tau & \\ \\
&= \ds  -\xi 
\frac{\partial}{\partial x} \int_0^{x}E_{\alpha}(-\lambda_0(x-\tau)^{\alpha}){\cal P}(\tau,t){\rm d}\tau 
 &\\ \\  
& = \ds -\xi {\cal P}(x,t) + \xi \int_0^x \lambda_0(x-\tau)^{\alpha-1}E_{\alpha,\alpha}[-\lambda_0(x-\tau)^{\alpha}]{\cal P}(\tau,t){\rm d}\tau  & \\ \\ 
& =\ds -\xi {\cal P}(x,t) + \xi \lambda_0 \mathfrak{E}_x (\alpha,\alpha; 1; -\lambda_0){\cal P}(x,t)\\ \\
\ds {\cal P}(x,t)|_{t=0} &= \delta(x) & 
\end{array}
\end{equation}
where in the second line it is used that
$-D_xE_{\alpha}(-\lambda_0x^{\alpha})=\lambda_0x^{\alpha-1}E_{\alpha,\alpha}(-\lambda_0x^{\alpha})$ which yields the (weakly singular and hence integrable) Mittag-Leffler density, and where $\mathfrak{E}_x$ is the Prabhakar integral acting on the space variable $x$ (see \cite{Prabhakar1971}). The integration limits reflect the upper triangular circulant property of
the Laplacian matrix function $g_{ML,\alpha}({\cal L})$ with
${\cal G}_{\xi,\lambda_0,\alpha}(x),{\cal W}_{ML,\alpha}(x,\lambda_0),{\cal P}(x,t)=0$ for $x<0$ (See also the last line of (\ref{Space-time-fractional-Bernoulli})). 
Equation (\ref{Cauchycont}) is the well-scaled continuous-space limit of (\ref{Space-time-fractional-Bernoulli}) and a
{\it generalized Kolmogorov-Feller (forward) diffusion limit equation} of general type (\ref{general}) solved by the state density kernel (\ref{finalresultspacecontilimit}). 
\\[1mm]
Eq. (\ref{Cauchycont}) has the following physical interpretation:
The Mittag-Leffler convolution on the right-hand side 
is the contribution to ${\cal P}(x,t)$ by {\it incoming jumps to $x$} which originate from 
all states $\tau$ with $0\leq \tau < x$. The term
$-\xi {\cal P}(x,t) = 
-\xi {\cal P}(x,t)\int_x^{\infty} {\cal W}_{ML,\alpha}(\tau-x){\rm d}\tau $ describes the ``loss'' due to {\it outgoing jumps from $x$} by long-range Mittag-Leffler jumps into the infinite half space $\tau>x$.
\\[1mm]
Coming back to Figure \ref{Fig2}: In view of the structure of Eq. (\ref{Cauchycont}), the emergence of oscillating behavior in the time dependence of the state density, especially
visible in the lower plot of this figure, is resulting from the complex interplay of incoming 
and outgoing jumps to and from a state $x$ where we give a rough qualitative interpretation: 
For $t$ and $x$ both `small', the incoming jumps at $x$ are `mainly' originating from the initial $\delta$-peak ${\cal P}(x,0)=\delta(x)$
where their accumulation giving rise to the 
first maximum which decays in time by outgoing jumps and the lack of incoming jumps. 
This effect becomes the more evanescent the larger $x$ due to
dispersion effects caused by jumps of any length drawn from the fat-tailed Mittag-Leffler density.
This dispersive behavior is strongly contrasted by the stable dispersion free propagation of the $\delta$-distribution state density in the standard Poisson process (see relation (\ref{specialcase11})).
\\[1mm]
It appears worthy also to consider briefly the connection 
with the Montroll-Weiss CTRW picture \cite{MontrollWeiss1965}. In this picture the Cauchy problem
(\ref{Space-time-fractional-Bernoulli}) is equivalent to a strictly increasing walk on the integer line
with discrete Mittag-Leffler jumps according to the one-step transition matrix (\ref{fracbernoulli-transmat})
subordinated to a fractional Poisson process with Mittag-Leffler waiting time density
(having time-Laplace transform 
${\tilde \chi}_{\beta}(s)=\frac{\xi}{\xi+s^{\beta}}$).
The generating function of the time-Laplace transform of the Cox series (\ref{coxseries}) then yields
the Montroll-Weiss equation 
\begin{equation}
\label{Montroll-Weiss}
\begin{array}{clc}
\ds {\bar p}_{\lambda,\xi}^{\alpha,\beta}(u,s) &= \ds \frac{1}{s}(1-{\tilde \chi}_{\beta}(s))\sum_{n=0}^{\infty} 
[{\tilde \chi}_{\beta}(s){\bar W}_{ML,\alpha}(u,\lambda)]^n & \\ \\  & =  \ds \frac{1}{s}(1-{\tilde \chi}_{\beta}(s))
\frac{1}{1-{\bar W}_{ML,\alpha}(u,\lambda){\tilde \chi}_{\beta}(s)} & \\ \\
& = \ds \frac{s^{\beta-1}}{s^{\beta}+ \frac{g_{ML,\alpha}(1-u,\lambda)}{g_{ML,\alpha}(1,\lambda)}}=  \ds \frac{s^{\beta-1}}{s^{\beta}+\frac{\xi(\lambda+1)(1-u)^{\alpha}}{\lambda+(1-u)^{\alpha}}} &
\end{array}
\end{equation}
where $g_{ML,\alpha}(1-u,\lambda)$ is the Laplacian generating function (\ref{generfulaplacmittagleff}) and ${\bar W}_{ML,\alpha}(u,\lambda)$ the generating function of the Mittag-Leffler transition matrix
(\ref{fracber}). We identify (\ref{Montroll-Weiss}) indeed with the time-Laplace transform of the Mittag-Leffler generating function  (\ref{generfufracber}) of the state-probabilities.
It is then straight-forward to see that (\ref{Montroll-Weiss}) is equivalent
to (\ref{writethegenfunction}) by rewriting Montroll-Weiss relation (\ref{Montroll-Weiss}) in the form
\begin{equation}
\label{rewritingintheform}
s^{\beta} {\bar p}(u,s) -s^{\beta-1}{\bar p}(u,t)|_{t=0}=  - \frac{\xi(\lambda+1)(1-u)^{\alpha}}{\lambda+(1-u)^{\alpha}} {\bar p}(u,s),\qquad {\bar p}(u,t)|_{t=0}  = 1.
\end{equation}
Inverting the time-Laplace transform on the left hand side yields the Caputo-derivative
thus the causal time-domain representation of (\ref{rewritingintheform}) indeed coincides with
the generating function representation of the Cauchy problem (\ref{writethegenfunction}).
This concludes our proof of equivalence of the space-time Mittag-Leffler process with a Montroll-Weiss CTRW of
a strictly increasing walk with discrete Mittag-Leffler jumps (\ref{fracbernoulli-transmat}) subordinated to a fractional Poisson process.
\section{GENERALIZED SPACE-TIME MITTAG-LEFFLER PROCESS}
\label{Prabhakargeneralisation}
Finally quasi as a byproduct we consider the space-fractional generalization of the space-time Mittag-Leffler process.
Let $g_{1,2}({\cal L})$ be good Laplacian matrix functions
constructed by L\'evy measures $\nu_{1,2}({\rm d}\tau)$ with (\ref{LKbernstein}). 
Then a good Laplacian  matrix function is obtained also by the chain $g_2(g_1({\cal L}))$. 
With the choice
$g_1({\cal L})=g_{ML,\alpha}({\cal L})$ and $g_2({\cal L})={\cal L}^{\mu}$ where ${\cal L}$ is Laplacian matrix (\ref{Laplacian}), we get ($\nu_2({\rm d}\tau)= -\frac{\tau^{-1-\mu}}{\Gamma(-\mu)}{\rm d}\tau$) 
\begin{equation}
\label{Prabhagen}
\begin{array}{clc}
\ds [g_{ML,\alpha}({\cal L})]^{\mu} &= \ds -\frac{1}{\Gamma(-\mu)}
\int_0^{\infty}(1-e^{-\tau \, g_{ML,\alpha}({\cal L})}) \tau^{-\mu-1}{\rm d}\tau ,\hspace{1cm} \mu \in (0,1) & \\ \\     
& = \ds  \frac{{\cal L}^{\alpha\mu}}{(\lambda {\mathbf 1}+{\cal L}^{\alpha})^{\mu}}. &
\end{array}
\end{equation}
This integral converges for $\mu \in (0,1)$ thus the fractional power $[g_{ML,\alpha}({\cal L})]^{\mu}$ in this range is a good Laplacian Bernstein matrix function retaining the Laplacian properties {\bf (i)-(iii)}. The Cauchy problem governing the state probabilities then reads
\begin{equation}
\label{Space-time-fractional-Bernoulli-Prabhakar-generalisation}
\begin{array}{lcl}
 \ds
\frac{d^{\beta}}{dt^{\beta}}p_n(t) &= &\ds   - \xi(\lambda+1)^{\mu}\frac{(1-{\hat T}_{-1})^{\alpha\mu}}{(\lambda+(1-{\hat T}_{-1})^{\alpha})^{\mu}}p_n(t)  ,\hspace{1cm} \xi,\, \lambda \, >\, 0, \hspace{0.25cm} t \geq 0 \hspace{0.25cm} \alpha,\, \mu, \, \beta \in (0,1] , \hspace{0.25cm} n \in \mathbb{N}_0
\\ \\  \ds p_n(t)|_{t=0} & = & \ds \delta_{0,n}
\end{array} 
\end{equation}
with the generating function of the state probabilities
\begin{equation}
\label{generfufracber-genprab}
{\bar p}^{\alpha,\mu,\beta}_{\lambda,\xi}(u,t) = \ds
E_{\beta}\left(\frac{-\xi t^{\beta} (\lambda+1)^{\mu}(1-u)^{\alpha\mu} }{(\lambda+(1-u)^{\alpha})^{\mu}}\right) = \ds \sum_{m=0}^{\infty}
 \frac{[-\xi t^{\beta} (\lambda+1)^{\mu}]^m}{\Gamma(\beta m+1)}[g_{ML,\alpha}(1-u)]^{m \mu}.
 \end{equation}
We hence get for the state-probabilities
\begin{equation}
\label{state-vector-Prabgen}
 p^{\alpha,\mu,\beta}_{\lambda,\xi,n}(t)  = \sum_{m=0}^{\infty}
 \frac{[- (\lambda+1)^{\mu}\xi t^{\beta}]^m}{\Gamma(\beta m+1)} 
 {\cal E}_{\alpha}^{(m\mu)}(\lambda,n) , \hspace{1cm} n \in \mathbb{N}_0
\end{equation}
with the cases (See also (\ref{prabhakarapprox}))
\begin{equation}\label{prabhakarapprox-gen}
 {\cal E}_{\alpha}^{(m\mu)}(\lambda,n) = \left\{\begin{array}{clc} \ds 
 \frac{1}{n!} \sum_{s=0}^{\infty} \frac{(-\lambda)^s (m\mu)_s}{s!} \frac{\Gamma(\alpha s+n)}{\Gamma(\alpha s)} , & 0<\lambda <1\\ \\ 
 \ds   \frac{(-1)^n\lambda^{-m\mu}}{n!} \sum_{s=0}^{\infty} \frac{(-1)^s \lambda^{-s} (m\mu)_s}{s!} \frac{\Gamma(\alpha(s + m\mu)+1)}{\Gamma(\alpha (s+m\mu)-n+1)} ,&  \lambda > 1
 \end{array} \right.
 \hspace{1cm} n \in \mathbb{N}_0
\end{equation}
where these series converge absolutely in view of the asymptotic behavior of the terms for large $s$, namely
$\sim s^{\mu m+n-1} \lambda^s$ ($\lambda <1$) and $\sim s^{\mu m+n-1} \lambda^{-s}$ ($\lambda > 1$).
${\cal E}_{\alpha}^{(m\mu)}(\lambda,n)$ is a discrete approximation of a Prabhakar kernel.
The well-scaled continuous-time limit
is now straight-forwardly obtained with $\lambda(h)=\lambda_0h^{\alpha} \rightarrow 0$ in the first line of (\ref{prabhakarapprox-gen}) with $n=\frac{x}{h}$ and yields the state density
\begin{equation}
\label{finalresultspacecontilimit-Prabgen}
\begin{array}{clc} \ds 
{\cal P}_{\lambda_0,\xi}^{\alpha,\beta}(x,t) & = \ds \lim_{h\to 0} \frac{1}{h} p^{\alpha,\mu,\beta}_{\lambda_0h^{\alpha},\xi,\frac{x}{h}}(t) \\ \\ & =  \ds  
\sum_{m=0}^{\infty}
 \frac{(-\xi t^{\beta})^m}{\Gamma(\beta m+1)}  e_{\alpha,0}^{m\mu}(-\lambda_0,x), \hspace{0.5cm} \lambda_0,\xi >0 ,\hspace{0.25cm} \alpha, \mu, \beta \in (0,1], \hspace{0.25cm} x,t \geq 0 & \\ \\
 &= \ds  D_x\left[\Theta(x) \Theta(t) \sum_{m=0}^{\infty}
 \frac{(-\xi t^{\beta})^m}{\Gamma(\beta m+1)} E_{\alpha,1}^{m\mu}(-\lambda_0x) \right] &
 \end{array}
\end{equation}
with the Prabhakar kernel $e_{\alpha,0}^{m\mu}(-\lambda_0,x)$ obtained as (See (\ref{toevaluate}))
\begin{equation}
\label{Prabker}
 e_{\alpha,0}^{m\mu}(-\lambda_0,x)  =  \ds \delta(x) + 
 \sum_{s=1}^{\infty}  \frac{ (m\mu)_s (-\lambda_0)^sx^{\alpha s-1} }{s!\Gamma(\alpha s)} , \hspace{1cm} m=0,1,2,.. \in \mathbb{N}_0 , \hspace{0.5cm} x \in  \mathbb{R}^{+}.
\end{equation}
Then we further get for the Laplacian density in the same way as in (\ref{space-contilim})
the kernel
\begin{equation}
\label{space-contilim-Prabgen}
{\cal G}_{\xi,\lambda_0,\alpha,\mu}(x) =  e_{\alpha,0}^{\mu}(-\lambda_0,x) =  \delta(x) - {\cal W}_{\alpha}^{(\mu)}(x) = D_x[\Theta(x) E_{\alpha,1}^{\mu}(-\lambda_0x^{\alpha})] 
\end{equation}
with $ E_{\alpha,1}^{\mu}(-\lambda_0x^{\alpha}) = \sum_{s=0}^{\infty}  \frac{ (\mu)_s (-\lambda_0)^sx^{\alpha s} }{s!\Gamma(\alpha s+1)}$. The transition density kernel hence yields
\begin{equation}
\label{kernelaxiliary}
{\cal W}_{\alpha}^{(\mu)}(x) = -D_xE_{\alpha,1}^{\mu}(-\lambda_0x^{\alpha}) = 
\lambda_0x^{\alpha-1}\sum_{s=0}^{\infty}\frac{(\mu)_{s+1}}{(s+1)!}
\frac{(-\lambda_0x^{\alpha})^s}{\Gamma(\alpha s +\alpha)} ,\hspace{1cm} x>0
\end{equation}
recovering for $\mu=1$ the Mittag-Leffler density. It is easily verified that the transition kernel is normalized $\int_0^{\infty}{\cal W}_{\alpha}^{(\mu)}(x){\rm d}x= E_{\alpha,1}^{\mu}(-\lambda_0x^{\alpha})|_{x=0}=1$.
The diffusion-limit Cauchy problem which again is of general type (\ref{general}) then reads

\begin{equation}
\label{contispacecauchy}
 \frac{\partial^{\beta}}{\partial t^{\beta}}{\cal P}(x,t) =-\xi {\cal P}(x,t) + 
\xi \int_0^x{\cal W}_{\alpha}^{(\mu)}(x-\tau){\cal P}(\tau,t){\rm d}\tau,\qquad {\cal P}(x,t)|_{t=0} = \ds \delta(x) 
\end{equation}
and is solved by the state density kernel (\ref{finalresultspacecontilimit-Prabgen}).
For $\mu=1$ all expressions turn into those previously derived for the space-time Mittag-Leffler process.
\\[1mm]
The fractional generalization of the space-time Mittag-Leffler process shows that the analysis of processes generated by chains
of Bernstein matrix functions $g_n(..g_1({\cal L}))$ which again are of the class of Bernstein functions retaining the good Laplacian properties may open an interesting direction to be further explored.
\section{CONCLUSIONS}
In this paper, we have analyzed space-time generalizations of the Poisson process defined by Cauchy problems
with generalized Kolmogorov-Feller difference-differential equations of type (\ref{Cauchyproblemgen}). 
These generalizations in the Montroll-Weiss CTRW picture are strictly 
increasing walks on the integer line time-changed with an independent renewal process. 
We have shown that this approach can also be applied more generally to biased walks on digraphs and
requires the construction of non-trivial `good Laplacian matrix functions' $g({\cal L})$. 
This introduces new topologies of fully connected structures with the small world property if the Laplacian matrix ${\cal L}$ is ergodic.
\\[1mm]
Choosing L\'evy densities related to normalized continuous distributions and a Laplacian ${\cal L}$ of the trivial strictly increasing walk (\ref{Laplacian})
leads to non-trivial
discrete-approximations of these L\'evy densities (see (\ref{LKbernstein})-(\ref{generfultransmat})) in the form of
{\it upper triangular circulant transition matrices} with strictly positive elements above the main diagonal, thus allowing any positive 
integer jump.
These properties are especially useful to construct new space-time generalizations of the Poisson
process.
\\[1mm]
As a pertinent example we have derived in this way a `good Laplacian matrix function' 
which generates a strictly increasing walk with discrete Mittag-Leffler jumps and introduced the
{\it space-time Mittag-Leffler process}. For this process, by means of explicit
formulae, we derived the state probabilities (Eq. (\ref{state-vector})) solving the Cauchy problem (\ref{Space-time-fractional-Bernoulli}).
Further, we developed a well-scaled continuous space limiting procedure and obtained the state density (\ref{finalresultspacecontilimit}) as a limiting expression of the state probabilities where Prabhakar kernels come into play.
We derived then the forward diffusion equation of general fractional type (\ref{Cauchycont}) which involves the Prabhakar integral and is solved by the state density (\ref{finalresultspacecontilimit}). The continuous-space limit diffusion equation (\ref{Cauchycont}) refers
to the general class of generalized Kolmogorov-Feller forward equations (\ref{general}). 
We also showed that the {\it space-time Mittag-Leffler process} in the 
Montroll-Weiss picture is a CTRW with discrete Mittag-Leffler jumps subordinated to an independent fractional Poisson process. 
\\[1mm]
Following this line, we introduced a space-fractional generalization of the space-time Mittag-Leffler process
constructed 
by a chain of two Laplacian Bernstein functions retaining the good Laplacian properties. In this way,
we derived the Cauchy problem (\ref{Space-time-fractional-Bernoulli-Prabhakar-generalisation}) governing the state probabilities for which we obtained an explicit formula (Eq. (\ref{state-vector-Prabgen})). We also deduced the well-scaled state density kernel (Eq. (\ref{finalresultspacecontilimit-Prabgen})) which involves Prabhakar kernels and 
solves the continuous-space Cauchy problem (\ref{contispacecauchy}).
\\[1mm]
There is a large potential in the presented approach of constructing space-time generalizations of the Poisson process in terms of strictly increasing walks. Also, new general types of {\it biased walks} time-changed with 
continuous-time or discrete-time counting processes derived in this way may have interesting applications in `birth-and-death' models.
On the other hand construction of new processes involving Prabhakar distributions (see e.g. \cite{PolitoScalas2016}) 
seem to be promising candidates due to their connections to the dynamics of certain complex phenomena. Further applications in
the field of `aging of complex systems' which are characterized
by strictly increasing random accumulation of damage (`misrepair') measures 
(see \cite{RiascosWangMi-Mi2029} for a Markovian model) could open an interesting
interdisciplinary field as well.

\subsection*{Acknowledgments}

F.~Polito has been partially supported by the project ``Memory in Evolving Graphs'' (Compagnia di San Paolo/Università degli Studi di Torino).
\begin{appendix}
\section{Discrete \texorpdfstring{$\bm{\delta}$}{delta}-distribution}
\label{Appendix}
We first introduce the {\it discrete 
Heaviside function} 
\begin{equation}
\label{heavisideintegers}
\ds
\Theta_h(x) = \Theta(x)= \left\{\begin{array}{l} 1 , \hspace{1cm} x\in \{ 0,h,2h,\ldots \hspace{0.2cm}\}\\ \\ 
                    0 ,\hspace{1cm}    x \in \{-h,-2h,\ldots \}
                   \end{array}\right. \hspace{1cm} x \in h\mathbb{Z}, \hspace{1cm} h>0.
                   \end{equation}
We especially emphasize that with this definition $\Theta_h(0)=1$.
We may extend the discrete Heaviside function here to its continuous counterpart defined for $x \in \mathbb{R}$,
i.e. the conventional Heaviside- unit-step function with 
$\Theta(x)=1$ for $t\geq 0$ (especially $\Theta(0)=1$) and $\Theta(x)=0$
for $x<0$. This is especially necessary when we use
${\hat T}_{-h} = e^{-hD_x}$ leading to the `distributional representation' 
$e^{-hD_x}\Theta(x)=\Theta(x-h)$ which is defined for $x\in \mathbb{R}$.
\\[1ex]
Let $\delta_{k,l}$ (we also use notation $\delta_{kl}$) be the circulant Kronecker-$\delta$ defined by
\begin{equation}
\label{Kroneckerdelta}
\ds 
\delta_{i,j} = \delta_{i+s,j+s}  \left\{\begin{array}{l} 1 , \hspace{1cm} i=j  \\ \\ 
                    0,     \hspace{1cm}   i \neq j
                   \end{array}\right. \hspace{2cm} i,j,s \in \mathbb{Z}.
                   \end{equation}
Then we define the `{\it discrete $\delta$-distribution}' as follows
\begin{equation}
\label{discrete-time-delta}
\ds \delta_h(x) = \frac{\Theta(x)-\Theta(x-h)}{h} = 
\frac{1-{\hat T}_{-h}}{h}\Theta(x) = \,\, \frac{1}{h}
\delta_{\frac{x}{h},0} \,\,  = \ds \left\{\begin{array}{l}\ds  \frac{1}{h} ,\hspace{0.5cm} x=0 \\ \\ \ds
0 ,\hspace{0.6cm} x\neq 0  \end{array}\right.
\hspace{1cm} x \in h\mathbb{Z} 
\end{equation}
where $\delta_{\frac{x}{h},0}$ in this relation indicates the circulant Kronecker-Symbol 
(\ref{Kroneckerdelta}) and $\frac{x}{h} \in \mathbb{Z}$. We observe that for $h=1$ we have $\delta_1(t)=\delta_{0,x}$.
Formula (\ref{discrete-time-delta}) defines a discrete density for
$x\in h\mathbb{Z}$. However, it makes sense to extend this definition to $x \in \mathbb{R}$. With this extended definition $\delta_h(x) = \frac{\Theta(x)-\Theta(x-h)}{h}$ becomes an integrable distribution (in the Gelfand-Shilov sense \cite{GelfandShilov1968}) with $\delta_h(x)= \frac{1}{h}$ for $x\in [0,h)$ and $\delta_h(x)=0$ else, especially we have $\int_{0}^{\infty}\delta_h(\tau){\rm d}\tau=\int_0^h\delta_h(\tau){\rm d}\tau=1$ 
thus
\begin{equation}
\label{dirac-delta}
\lim_{h\rightarrow 0} \delta_h(x) = 
\lim_{h\rightarrow 0}\frac{1-e^{-hD_x}}{h} \Theta(x) = D_x\Theta(x) = \delta(x) ,\hspace{1cm} x \in \mathbb{R}
\end{equation}
is a non-symmetric Dirac's $\delta$-distribution which is concentrated at $0+$ (and is null at $0-$), fulfilling therefore $\int_0^{\infty}\delta(x){\rm d}x = 
\lim_{h\to 0} \int_{0}^{h}\delta_h(x){\rm d}x =1$. This property is absolutely crucial and ensures for instance the normalization of the state density kernel (\ref{finalresultspacecontilimit}).
It is also worthy to consider the (spatial) Laplace transform
\begin{equation}
\label{Lapladeltah}
\begin{array}{l} \ds
{\cal L}\{\delta_h(x)\}(s) = \int_{0_{-}}^{\infty} \delta_h(x) e^{-sx}{\rm d}x = \int_{-\infty}^{\infty} e^{-sx} \frac{(1-e^{-hD_x})}{h}\Theta(x){\rm d}x \\ \\ \ds
\hspace{1cm} = \int_{-\infty}^{\infty} \Theta(x) \frac{(1-e^{+hD_x})}{h} e^{-sx} {\rm d}x =\frac{1-e^{-hs}}{h} \int_0^{\infty}e^{-sx} {\rm d}x \hspace{1cm} \Re\{s\} >0\\ \\ \ds \hspace{1cm}
 = \frac{1-e^{-hs}}{h s} 
\end{array}
\end{equation}
where indeed $\lim_{h\rightarrow 0} {\cal L}\{\delta_h(x)\}(s) = {\cal L}\{\delta(x)\}(s) = 1$ 
recovers the Laplace transform of Dirac's $\delta$-distribution (\ref{dirac-delta}).
\end{appendix}

\end{document}